%% file: SimplSuperrigid.tex
\documentclass[11pt]{amsart}
\usepackage{latexsym, amssymb, amscd, mathrsfs, graphics, fullpage}
\usepackage[backref=page,colorlinks,bookmarksopen=true]{hyperref}

\newtheorem{theorem}{Theorem} 
\newtheorem{cor}[theorem]{Corollary}
\newtheorem{lem}[theorem]{Lemma}
\newtheorem{prop}[theorem]{Proposition}
\newtheorem{defi}[theorem]{Definition}

\theoremstyle{definition}
\newtheorem*{remark}{Remark}
\newtheorem*{remarks}{Remarks}

\newcommand{\eps}{\epsilon}
\newcommand{\la}{\langle}
\newcommand{\ra}{\rangle}
\newcommand{\inv}{^{-1}}
\newcommand{\E}{\mathbb{E}}
\newcommand{\F}{\mathbf{F}}
\newcommand{\K}{\mathbf{K}}
\newcommand{\N}{\mathbf{N}}
\newcommand{\R}{\mathbf{R}}
\newcommand{\Z}{\mathbf{Z}}
\newcommand{\GL}{\mathrm{GL}}
\newcommand{\SL}{\mathrm{SL}}

\DeclareMathOperator{\Stab}{Stab}  \DeclareMathOperator{\Ker}{Ker}
\DeclareMathOperator{\Hom}{Hom} \DeclareMathOperator{\Aut}{Aut} \DeclareMathOperator{\Isom}{Isom}
\DeclareMathOperator{\Pc}{Pc} \DeclareMathOperator{\Out}{Out} \DeclareMathOperator{\Inn}{Inn}
\DeclareMathOperator{\proj}{proj}

\def\min{\mathop{\mathrm{min}}\nolimits}
\def\max{\mathop{\mathrm{max}}\nolimits}

\def\og{\leavevmode\raise.3ex\hbox{$\scriptscriptstyle\langle\! \langle\,$}}
\def\fg{\leavevmode\raise.3ex\hbox{$\scriptscriptstyle\,\rangle\! \rangle\ $}}

\begin{document}

\title[Simplicity and superrigidity of twin building lattices]
{Simplicity  and superrigidity of twin building lattices}
\author{Pierre-Emmanuel Caprace$^*$ and Bertrand R\'emy}\thanks{$^*$ F.N.R.S. research fellow}
\date{\today}
\maketitle

\vspace{1cm}

{\footnotesize {\bf Abstract:}~ Kac-Moody groups over finite fields are finitely generated groups. Most of them
can naturally be viewed as irreducible lattices in products of two closed automorphism groups of non-positively
curved twinned buildings: those are the most important (but not the only) examples of \emph{twin building
lattices}. We prove that these lattices are simple if and only if the corresponding buildings are (irreducible
and) not of affine type (i.e. they are not Bruhat-Tits buildings). In fact, many of them are finitely presented
and enjoy property (T). Our arguments explain geometrically why simplicity fails to hold only for affine
Kac-Moody groups. Moreover we prove that a nontrivial continuous homomorphism from a completed Kac-Moody group
is always proper. We also show that Kac-Moody lattices fulfill conditions implying strong superrigidity
properties for isometric actions on non-positively curved metric spaces. Most results apply to the general class
of twin building lattices.

\smallskip


}

\vspace{5mmplus1mmminus1mm}

\section*{Introduction}
\label{s - intro}

Since the origin, Kac-Moody groups (both in their so-called minimal and maximal versions) have been mostly
considered as natural analogues of semisimple algebraic groups arising in an infinite-dimensional Lie theoretic
context (see e.g. \cite{KP83} and \cite{KumarBook}). A good illustration of this analogy is the construction of
minimal Kac-Moody groups over arbitrary fields, due to J.~Tits \cite{TitsJA}, by means of presentations
generalizing in infinite dimension the so-called Steinberg presentations of Chevalley groups over fields
\cite{Steinberg}. This presentation provides not only a Kac-Moody group $G$, but also a family of root subgroups
$(U_\alpha)_{\alpha \in \Phi}$ indexed by an abstract root system $\Phi$ and satisfying a list of properties
shared by the system of root groups of any isotropic semisimple algebraic group. These properties constitute the
group theoretic counterpart of the geometric notion of a twin building: any group endowed with such a family of
root groups, which is called a \emph{twin root datum}, has a natural diagonal action on a product of two
buildings, and this action preserves a twinning. We refer to \cite{TitsLMS} and \cite{RemAst} for this
combinatorial point of view.

\vspace{2mm}

In this paper, we are mainly interested in finitely generated Kac-Moody groups, i.e. minimal Kac-Moody groups
over finite ground fields. In this special situation, it has been noticed more recently that another viewpoint,
different from the aforementioned algebraic group theoretic one, is especially relevant: the arithmetic group
viewpoint. The striking feature which justifies in the first place this more recent analogy is the fact that
finitely generated Kac-Moody groups embed as irreducible lattices in the product of closed automorphism groups
of the associated buildings, provided the ground field is sufficiently large, see \cite{RemCRAS}. A sufficient
condition for this is that the order of the finite ground field is at least the number of canonical generators
of the Weyl group. In fact, Kac-Moody theory is one of the few known sources of examples of irreducible lattices
in products of locally compact groups outside the classical world of lattices in higher-rank Lie groups. On the
other hand, the intersection between Kac-Moody groups and arithmetic groups is nonempty since Kac-Moody groups
of affine type, namely those obtained by evaluating Chevalley group schemes over rings of Laurent polynomials,
are indeed arithmetic groups. A standard example, which is good to keep in mind, is the arithmetic group
$\SL_n(\F_q[t, t\inv])$, which is an irreducible lattice of $\SL_n \bigl(\F_q (\!(t)\!)\bigr) \times \SL_n
\bigl(\F_q (\!(t\inv)\!)\bigr)$. This arithmetic group analogy, suggesting the existence of strong similarities
between Kac-Moody groups of arbitrary type and the previous examples of affine type, is supported by several
other results, see e.g. \cite{Ab95} for finiteness properties, \cite{DJ02} for continuous cohomology,
\cite[\S1]{RR06} for some structural properties, etc.

\vspace{2mm}

The main result of the present paper shows that for infinite Kac-Moody groups over finite fields, there is a
sharp structural contrast between  affine and non-affine groups. Indeed, affine Kac-Moody groups over finite
fields, as finitely generated linear groups, are residually finite. On the other hand, non-affine Kac-Moody
groups are submitted to the following:

\vspace{2mm}

{\bf Simplicity theorem (Kac-Moody version).}\it~
Let $\Lambda$ be a split or almost split Kac-Moody group over
a finite  field ${\bf F}_q$. Assume that the Weyl group $W$ of $\Lambda$ is an irreducible, infinite and
non-affine Coxeter group. Then every finite index subgroup of $\Lambda$ contains the derived subgroup $[\Lambda,
\Lambda]$, which is of finite index. Assume moreover that $q \geq |S|$. Then the group $[\Lambda, \Lambda]$,
divided by its finite center, is simple. \rm

\vspace{2mm}

A more general result (Theorem \ref{thm:SimpleDRJ} of Subsect. \ref{ss - simple lattices}) holds in the abstract
framework of twin root data; it was announced in \cite{CapRem}.

\vspace{2mm}

It follows from the above that for any neither spherical nor affine, indecomposable  generalized Cartan matrix
$A$ of size $n$, there exists a Kac-Moody group functor  $\mathcal{G}_A$ such that for any finite field ${\bf
F}_q$ with $q \geq n > 2$, the group $\Lambda=\mathcal{G}_A({\bf  F}_q)$, divided by its finite center, is an
infinite finitely generated simple group. We also note that this simplicity result for Kac-Moody groups over
finite fields implies strong non-linearity properties for Kac-Moody groups over \emph{arbitrary} fields of
positive characteristic (Theorem \ref{thm:NonLinearity}). To be more constructive, we add the following
corollary (see Corollary \ref{cor - T and simple}). As pointed out to us by Y.~Shalom, we obtain the first
infinite {\it finitely presented}~discrete groups to be both simple and Kazhdan.

\vspace{2mm}

{\bf Simple Kazhdzan group corollary.}\it~ If the generalized Cartan matrix $A$ is $2$--spherical (i.e. every $2
\times 2$--submatrix is of spherical type) and if $q > 1764^n$, then the group $\Lambda/Z(\Lambda)$ is finitely
presented, simple and Kazhdan. Moreover there exist infinitely many isomorphism classes of infinite groups with
these three properties. \rm

\vspace{2mm}

Another consequence is the possibility to exhibit a large family of inclusions of lattices in topological groups
for which the density of the commensurator does not hold (see Corollary \ref{cor - non arithmetic}).

\vspace{2mm}

{\bf Non-arithmeticity corollary.}\it~ Let $\Lambda$ be a split or almost split Kac-Moody group over a finite
field ${\bf F}_q$. Assume that the Weyl group of $\Lambda$ is irreducible, infinite and non-affine. Let
$\mathscr{B}_+, \mathscr{B}_-$ be the buildings associated with $\Lambda$ and let $\overline\Lambda_+$ and
$\overline\Lambda_-$ be the respective closures of the images of the natural actions $\Lambda \to
\Aut(\mathscr{B}_+)$ and $\Lambda \to \Aut(\mathscr{B}_-)$. We view $\Lambda/Z(\Lambda)$ as a diagonally
embedded subgroup of $\overline\Lambda_-\times\overline\Lambda_+$. Then the commensurator ${\rm
Comm}_{\overline\Lambda_-\times\overline\Lambda_+}\bigl( \Lambda/Z(\Lambda) \bigr)$ contains $\Lambda$ as a
finite index subgroup; in particular it is discrete. \rm

\vspace{2mm}

The proof of the above theorem follows a two--step strategy which owes much to a general approach initiated by
M.~Burger and Sh.~Mozes \cite{BMProducts} to construct simple groups as cocompact lattices in products of two
trees; however, the details of arguments are often substantially different. The idea of  \cite{BMProducts} is to
prove the normal subgroup property (i.e. normal subgroups either are finite and central or have finite index)
following G.~A.~Margulis' proof for lattices in higher rank Lie groups \cite[VIII.2]{Margulis} but to disprove
residual finiteness by geometric arguments in suitable cases \cite[Proposition 2.1]{BMProducts}. This relies on
a preliminary study of sufficiently transitive groups of tree automorphisms \cite{BMTrees}. Let us also recall
that a finitely generated just infinite group (i.e. all of whose proper quotients are finite) is either
residually finite or is, up to finite index, a direct product of finitely many isomorphic simple groups
\cite{Wilson}.

\vspace{2mm}

In our situation, the first step of the proof, i.e. the normal subgroup property, had been established in
previous papers, in collaboration with U.~Bader and Y.~Shalom: \cite{BSh06}, \cite{Rem05} and \cite{Shalom}.
This fact, which is recalled here as Theorem \ref{th - NSP}, is one of the main results supporting the analogy
with arithmetic groups mentioned above. One difference with \cite[Theorem 4.1]{BMProducts} is that the proof
does not rely on the Howe--Moore property (i.e. decay of matrix coefficients). Instead, Y.~Shalom and U.~Bader
use cohomology with unitary coefficients and Poisson boundaries. In fact, it can be seen that closed strongly
transitive automorphism groups of buildings do not enjoy Howe--Moore property in general: whenever the ambient
closed automorphism group of the buildings in consideration contains a proper parabolic subgroup which is not of
spherical type (i.e. whose Weyl subgroup is infinite), then any such parabolic subgroup is an open subgroup
which is neither compact nor of finite index.

\vspace{2mm}

The second half of the simplicity proof does not actually use the notion of residual finiteness. Instead, it
establishes severe restrictions on the existence of finite quotients of a group endowed with a twin root datum
of non-affine type (Theorem \ref{thm:NoFiniteQuotient}). Here, one confronts the properties of the system of
root subgroups to a geometric criterion which distinguishes between the Tits cones of affine (i.e. virtually
abelian) and of non-affine Coxeter groups (see \cite{MarVin00} and \cite{NosVin02} for another illustration of
this fact). This part of the arguments holds without any restriction on the ground field, and holds in
particular for those groups over tiny fields for which simplicity is still an open question.

\vspace{2mm}

We note that the simplicity theorem above is thus obtained as the combination of two results which pertain
respectively to each of the two analogies mentioned above. In this respect, it seems that the structure of
Kac-Moody groups is enriched by the ambiguous nature of these groups, which are simultaneously
arithmetic-group-like and algebraic-group-like.

\vspace{2mm}

In order to conclude the presentation of our simplicity results, let us compare quickly Burger-Mozes' groups
with simple Kac-Moody lattices. The groups constructed thanks to \cite[Theorem 5.5]{BMProducts} are cocompact
lattices in a product of two trees; they are always finitely presented, simple, torsion free and amalgams of
free groups (hence cannot have property (T)). Simple Kac-Moody lattices are non-uniform lattices of products of
buildings, possibly (and usually) of dimension $\ge 2$; they are often (not always, though) finitely presented
and Kazhdan and contain infinite subgroups of finite exponent. It is still an open and challenging question to
construct simple cocompact lattices in higher-dimensional buildings.

\vspace{2mm}

We now turn to the second series of results of this paper. It deals with restrictions on actions of Kac-Moody
lattices on various non-positively curved spaces. This is a natural question with respect to the previous
results in view of  the following known fact: a non residually finite group cannot be embedded nontrivially into
a compact group (see Proposition \ref{prop - homomorphisms}). In particular a simple Kazhdan group acting
non-trivially on a Gromov-hyperbolic proper metric space $Y$ of bounded geometry cannot fix any point in the
visual compactification $\overline Y = Y \sqcup \partial_\infty Y$. This is where superrigidity results come
into play since they provide even stronger restrictions on actions: a finitely generated group of Kac-Moody type
is naturally a lattice for a product of two locally compact totally disconnected groups. The latter groups are
strongly transitive closed automorphism groups of buildings and as such have a $BN$-pair structure.

\vspace{2mm}

Therefore, it makes sense to try to use the recent superrigidity results due to N.~Monod \cite{Monod06} and to
T.~Gelander, A.~Karlsson and G.~Margulis \cite{GKM06}. As for Y. Shalom's result on property (T) \cite{Shalom},
the non-cocompactness of Kac-Moody lattices is an obstruction to a plain application of these results (stated
for uniform lattices). Still, all the previous references propose measure-theoretic or representation-theoretic
substitutes for cocompactness. Although weak cocompactness of Kac-Moody groups is still open, we prove a
substitute proposed in \cite{GKM06} (see Theorem \ref{th - p-integrable}).

\vspace{2mm}

{\bf Uniform integrability theorem.}\it~ Let $\Lambda$ be a split or almost split Kac-Moody group over $\F_q$.
Assume that $\Lambda$ is a lattice of the product of its twinned buildings $\mathscr{B}_\pm$. Then the group
$\Lambda$ admits a natural fundamental domain with respect to which it is uniformly $p$--integrable for any $p
\in [1;+\infty)$. \rm \vspace{2mm}

For an arbitrary inclusion of a finitely generated lattice $\Gamma$ in a locally compact group $G$, uniform
integrability is a technical condition requiring the existence of a fundamental domain $D$ with respect to which
some associated cocycle is uniformly integrable (see Subsect. \ref{ss - uniform integrability}); it ensures the
existence of harmonic maps \cite{GKM06}. Here is an example of an available superrigidity statement
\cite[Theorem 1.1]{GKM06}.

\vspace{2mm}

{\bf Superrigidity corollary.}\it~ Let $\Lambda$ be as above; in particular it is a lattice of the product of
its two completions $\overline\Lambda_- \times \overline\Lambda_+$. Let $X$ be a complete Busemann
non-positively curved uniformly convex metric space without nontrivial Clifford isometries. We assume that there
exists a $\Lambda$-action by isometries on $X$ with reduced unbounded image.
Then the $\Lambda$-action extends uniquely to a $\overline\Lambda_- \times \overline\Lambda_+$-action which
factors through $\overline\Lambda_-$ or $\overline\Lambda_+$. \rm\vspace{2mm}

The relevancy of reduced actions was pointed out in \cite{Monod06}: a subgroup $L<{\rm Isom}(X)$ is called
\textit{reduced} if there is no unbounded closed convex proper subset $Y$ of $X$ such that $gY$ is at finite
Hausdorff distance from $Y$ for any $g  \in  L$. We also recall that a \textit{Clifford isometry} of $X$ is a
surjective isometry $T : X \to X$ such that $x \mapsto d(T(x),x)$ is constant on $X$.

\vspace{2mm}

These results about continuous extensions of group homomorphisms call  for structure results for the geometric
completions $\overline\Lambda_\pm$ of  Kac-Moody groups over finite fields, i.e. the closures of the
non-discrete $\Lambda$-actions on each building $\mathscr{B}_\pm$. Indeed, once a continuous extension has been
obtained by superrigidity, it is highly desirable to determine whether  this map is proper, e.g. to know whether
infinite discrete subgroups can have a global fixed point in the target metric space. When the ambient
topological groups are semisimple Lie groups, the properness comes as a consequence of the Cartan decomposition
of such groups \cite[Lemma 5.3]{BM96}. The difficulty in our situation is that, with respect to structure
properties, topological groups of Kac-Moody type are not as nice as  semisimple algebraic groups over local
fields. Unless the Kac-Moody group is of affine type, there is no Cartan decomposition in which double cosets
modulo a maximal compact subgroup are indexed by an abelian semi-group: the Weyl group is not virtually abelian
and  roots cannot be put into finitely many subsets according to parallelism classes of walls in the Coxeter
complex. This is another avatar of the strong Tits alternative for infinite Coxeter groups \cite{MarVin00},
\cite{NosVin02}. Here is a slightly simplified version of our main properness result
(Theorem~\ref{th:proper}).

\vspace{2mm}

{\bf Properness theorem.}\it~
Let $\Lambda$ be a split or almost split simply connected Kac-Moody group over
$\F_q$ and let $\overline\Lambda_+$ be its positive topological completion. Then any nontrivial continuous
homomorphism $\varphi : \overline\Lambda_+ \to G$ to a locally compact second countable group $G$ is proper. \rm

\vspace{2mm}

As an example of application of superrigidity results, we study actions on Kac-Moody lattices on
$\mathrm{CAT}(-1)$-spaces. In this specific case, the most appropriate results available are the superrigidity
theorems of N.~Monod and Y.~Shalom \cite{MS04}. Putting these together with the abstract simplicity of
non-affine Kac-Moody lattices and the properness theorem above enables us to exhibit strong incompatibilities
between higher-rank Kac-Moody groups and some negatively curved metric spaces (see Theorem \ref{th - action on
CAT(-1)} for more details).

\vspace{2mm}

{\bf \og Higher-rank versus ${\bf {\rm CAT}(-1)}$\fg theorem.}\it~ Let $\Lambda$ be a simple Kac-Moody lattice
and $Y$ be a  proper ${\rm CAT(-1)}$-space with cocompact isometry group. If the buildings $\mathscr{B}_\pm$
of $\Lambda$ contain flat subspaces of dimension $\geq 2$ and if $\Lambda$ is Kazhdan, then the  group $\Lambda$
admits no nontrivial action by isometries on $Y$. \rm\vspace{2mm}

We show below, by means of a specific example, that the assumption that $\Isom(Y)$ is cocompact is necessary
(see the remark following Theorem~\ref{th - action on CAT(-1)}). This theorem was motivated by
\cite[Corollary~0.5]{BM96}. Note that we made two assumptions (one on flat subspaces, one on property~(T))
which, in the classical case, are implied by the same algebraic condition. Namely, if $\Lambda$ were an
irreducible lattice in a product of semisimple algebraic groups, and if each algebraic group were of split rank
$\ge 2$, then both \og higher-rank\fg assumptions would be fulfilled. In the Kac-Moody case, there is no
connection between existence of flats in the building and property (T). The relevant rank here is the geometric
one (the one involving flats in the buildings). According to \cite{BRW05} and \cite{CapraceHaglund}, it has a
more abstract interpretation relevant to the general theory of totally disconnected locally compact groups.

\vspace{2mm}

The proofs of most results of the present paper use in a very soft way that the construction of the lattices
considered in this introduction pertains to Kac-Moody theory. The actual tool which is the most natural
framework for our arguments is the notion of a \emph{twin root datum} introduced in~\cite{TitsLMS}. It turns out
that the class of groups endowed with a twin root datum includes split and almost split Kac-Moody groups only as
a (presumably small) subfamily (see \S\ref{s - KM lattices}  below). Several exotic constructions of such groups
outside the strict Kac-Moody framework are known, see e.g. \cite[\S 9]{Ti88/89} for groups acting on twin trees,
\cite{RR06} for groups acting on right-angled twin buildings and \cite{Mu99} for groups obtained by integration
of Moufang foundations. All these examples are discrete subgroups of the product of the automorphism groups of
the two halves of a twin building, which are actually mostly of finite covolume. These lattices, called
\emph{twin building lattices}, constitute the main object of study for the rest of this paper.

\vspace{2mm}

{\bf Structure of the paper.}~ In the preliminary Sect.~\ref{s - prel}, we fix the conventions and notation.
Sect.~\ref{s - KM lattices} is devoted to collect some basic results for later reference. Although these results
are often stated in the strict Kac-Moody framework in the literature, we have been careful to state and prove
them in the context of twin building lattices. Sect.~\ref{s - AptTopology} introduces a fixed-point property of
root subgroups and it is shown that most examples of twin building lattices enjoy this property. It is then used
to establish several useful structural properties of these completions. In Sect. \ref{s - non-affine Coxeter},
we prove the main fact needed for the simplicity theorem; it is the existence of a weakly hyperbolic geometric
configuration of walls for non-affine infinite Coxeter complexes. In Sect. \ref{s - simple lattices}, the
simplicity theorem is proved together with very strong restrictions on quotients of the groups for  which the
simplicity is still unknown. In Sect. \ref{s - non linear lattices}, we prove a non-linearity property for
Kac-Moody groups over arbitrary fields of positive characteristic. In Sect. \ref{s - homomorphisms}, we study
homomorphisms from Kac-Moody groups to arbitrary locally compact groups; the main part of this section deals
with the geometric completions of Kac-Moody groups, more  precisely it establishes that any nontrivial
continuous homomorphism on must be proper. In Sect. \ref{s - superrigidity}, we prove that Kac-Moody groups
enjoy  the uniform integrability condition required by Gelander-Karlsson-Margulis to derive superrigidity
statements; restrictions on actions on hyperbolic metric spaces in terms of \og rank\fg are derived  from this.

\vspace{2mm}

{\bf Acknowledgements.}~ We thank U.~Baumgartner, M.~Burger, M.~Ershov, T.~Gelander, N.~Monod, B.~M\"uhlherr,
Y.~Shalom for helpful discussions. The first author thanks the F.N.R.S. for supporting a visit to the
Universit{\'e} de Lyon~1, where most of this work was accomplished. The second author thanks the
Max-Planck-Institut f\"ur Mathematik in Bonn for its hospitality while a final version of this paper has been
written.

\vspace{1cm}

\setcounter{section}{-1}
\section{Notation and general references}\label{s - prel}
Let us fix some notation, conventions and make explicit our standard references.

\vspace{2mm}

\subsection{About Coxeter groups} Throughout this paper, $(W,S)$ denotes  a Coxeter system
\cite[IV.1]{BbkLie4-5-6} of finite rank  (i.e. with $S$ finite) and $\ell$ or $\ell_S$ denotes the word length
$W \to \N$ with respect to the generating set $S$. We denote by $W(t)$ the canonical growth series, i.e. the
series  $\sum_{n \geq 0} c_nt^n$ where $c_n$ is the number of elements  $w \in W$ such that $\ell_S(w)=n$. The
combinatorial root system $\Phi$ of $W$ is abstractly defined in \cite[Sect. 5]{TitsJA}. We adopt this point of
view because, since it is purely set-theoretic, it is useful to connect several geometric realizations of the
Coxeter complex of $(W,S)$ \cite{Brown}, \cite{Ronan}. A pair of opposite roots here is a pair of complementary
subsets $W$ which are permuted by a suitable conjugate of some canonical generator $s \in  S$. The set of simple
roots is denoted by $\Pi$.

\vspace{2mm}

Recall that a set of roots $\Psi$ is called \textit{prenilpotent} if both intersections $\bigcap_{\alpha\in\Psi}
\alpha$ and $\bigcap_{\alpha\in\Psi} -\alpha$ are nonempty. Given a prenilpotent pair $\{\alpha, \beta \}
\subset \Phi$, we set

\vspace{2mm}

\centerline{$[\alpha, \beta] := \{ \gamma \in \Phi \; | \; \alpha \cap \beta \subset \gamma \text{ and }
(-\alpha) \cap (-\beta) \subset -\gamma \}$.}

\subsection{About geometric realizations} We denote by $\mathscr{A}$ the Davis complex associated to $(W,S)$
and by $d$ the corresponding CAT(0) distance on $\mathscr{A}$ \cite{Davis}. The metric space $\mathscr{A}$ is
obtained as a gluing of compact subsets, all isometric to one another and called chambers. The group $W$ acts
properly discontinuously on $\mathscr{A}$ and simply transitively on the chambers. The fixed point set of each
reflection, i.e. of each element of the form $wsw\inv$ for some $s  \in  S$ and $w  \in  W$, separates
$\mathscr{A}$ into two disjoint halves, the closures of which are called root half-spaces. The fixed point set
of a reflection is called a wall. The set of root half-spaces of $\mathscr{A}$ is denoted by
$\Phi(\mathscr{A})$; it is naturally in $W$--equivariant bijection with $\Phi$. We distinguish a base chamber,
say $c_+$, which we call the standard chamber: it corresponds to $1_W$ in the above free $W$-action. We denote
by $\Phi_+$ the set of root half-spaces containing $c_+$ and by $\Pi$ the set of simple roots, i.e. of positive
roots bounded by a wall associated to some $s  \in  S$.

\subsection{About group combinatorics}\label{sect:TRD}
The natural abstract framework in which the main results of this paper hold is provided by the notion of a
\emph{twin root datum}, which was introduced in \cite{TitsLMS} and is further discussed for instance in \cite[\S
1]{Ab95} or \cite[1.5]{RemAst}. A twin root datum consists of a couple $(G,\{U_\alpha\}_{\alpha\in\Phi})$ where
$G$ is a group and $\{U_\alpha\}_{\alpha\in\Phi}$ is a collection of subgroups indexed by the combinatorial root
system of some Coxeter system; the subgroups $\{U_\alpha\}_{\alpha\in\Phi}$, called \emph{root groups}, are
submitted to the following axioms, where  $T := \bigcap_{\alpha\in\Phi} N_G(U_\alpha)$ and $U_+$ (resp. $U_-$)
denotes the subgroup generated by the root groups indexed by the positive roots (resp. their opposites):
\begin{description}
\item[(TRD0)] For each $\alpha \in \Phi$, we have $U_\alpha \neq \{1\}$.

\item[(TRD1)] For each prenilpotent pair $\{\alpha, \beta\} \subset \Phi$, the commutator group $[U_\alpha,
U_\beta]$ is contained in the group $U_{]\alpha, \beta[}:= \la U_\gamma | \; \gamma \in ]\alpha, \beta[ \ra$.

\item[(TRD2)] For each $\alpha \in \Pi$ and each $u \in U_\alpha \backslash \{1\}$, there exists elements $u',
u'' \in U_{-\alpha}$ such that the product $m(u):=u' u u''$ conjugates $U_\beta$ onto $U_{s_\alpha(\beta)}$ for
each $\beta \in \Phi$.

\item[(TRD3)] For each $\alpha \in \Pi$, the group $U_{-\alpha}$ is not contained in $U_+$ and the group
$U_{\alpha}$ is not contained in $U_-$.

\item[(TRD4)] $G = T \la U_\alpha |\; \alpha \in \Phi \ra$.
\end{description}

We also set $N := T . \la m(u) \; | \; u  \in  U_\alpha-\{1\}, \, \alpha \in \Pi \ra$. A basic fact is that the
subquotient $N/T$ is isomorphic to $W$; we call it the \emph{Weyl group} of $G$.

\subsection{About twin buildings} The geometric counterpart to twin root data is the notion of twin
buildings. Some references are \cite{TitsLMS},  \cite[\S 2]{Ab95} or \cite[\S 2.5]{RemAst}. Roughly speaking, a
group with a twin root datum $\{U_\alpha\}_{\alpha\in\Phi}$ of type $(W,S)$ admits two structures of $BN$--pairs
which are not conjugate to one another in general. Let $(\mathscr{B}_+, \mathscr{B}_-)$ be the associated
twinned buildings; their apartments are modelled on the Coxeter complex of $(W,S)$. We will not need the
combinatorial notion of a twinning between $\mathscr{B}_-$ and $\mathscr{B}_+$. The standard twin apartment
(resp. standard positive chamber) is denoted by $(\mathscr{A}_+, \mathscr{A}_-)$ (resp. $c_+$). We identify the
Davis complex $\mathscr{A}$ with the positive apartment  $\mathscr{A}_+$. With this identification and when the
root groups are all finite, the buildings $\mathscr{B}_\pm$ are locally finite ${\rm CAT}(0)$ cell complexes.

\vspace{1cm}

\section{Twin building lattices and their topological completions}\label{s - KM lattices}

As mentioned in the introduction, the main results of this paper apply not only to split or almost split
Kac-Moody groups over finite fields, but also to the larger class of groups endowed with a twin root datum with
finite root groups. Some existing results in the literature are stated for split or almost split Kac-Moody
groups, but remain actually valid in this more general context of twin root data. The purpose of this section is
to collect some of these results and to restate them in this context for subsequent references.

\vspace{2mm}

\subsection{Kac-Moody groups versus groups with a twin root datum}\label{ss - exoticTRD}
Although the notion of a twin root datum was initially designed as an appropriate tool to study Kac-Moody
groups, it became rapidly clear that many examples of twin root data arise beyond the strict scope of Kac-Moody
theory. This stands in sharp contrast to the finite-dimensional situation: as follows from the classification
achieved in \cite{Ti74} (see also \cite{TW02}), any group endowed with a twin root datum with \emph{finite} Weyl
group of rank at least~$3$ and of irreducible type is associated (in a way which we will not make precise) with
some isotropic simple algebraic group over a field or with a classical group over a (possibly skew) field. Here
is a list of known constructions which yields examples of twin root data with \emph{infinite} Weyl group but not
associated with split or almost split Kac-Moody groups:

\begin{itemize}
\item[(I)] \cite[\S 9]{Ti88/89} constructs a twin root datum with infinite dihedral Weyl group and arbitrary
prescribed rank one Levi factors. The possibility of mixing ground fields prevent these groups from being of \og
Kac-Moody origin\fg. The associated buildings are one-dimensional, i.e. trees.

\item[(II)] In \cite{RR06}, the previous construction is generalized to the case of Weyl groups which are
arbitrary right-angled Coxeter groups. In particular, the associated buildings are of arbitrarily large
dimension.

\item[(III)] Opposite to right-angled Coxeter groups are $2$--spherical Coxeter groups, i.e. those Coxeter
groups in which every pair of canonical generators generates a finite group. Twin root data with $2$--spherical
Weyl group are submitted to strong structural restrictions (see \cite{MR94}) showing in particular that wild
constructions as in the right-angled case are impossible. For instance, the following fact is a consequence of
the main result of \cite{Mu99}: \emph{a group $\Lambda$ endowed with a twin root datum $\{U_\alpha\}_{\alpha \in
\Phi}$ of irreducible type, such that the root groups are all finite of order~$> 3$ and generate $\Lambda$, and
that every rank~$2$ parabolic subgroup is of type $A_1 \times A_1$, $A_2$, $B_2$, $C_2$ or $G_2$ must be a split
or almost split Kac-Moody group in the sense of \cite{RemAst}}. Furthermore, it was mentioned to us by
B.~M\"uhlherr, as a non-obvious strengthening of \cite{Mu99}, that the preceding statement remains true if one
allows the rank~$2$ subgroups to be twisted Chevalley groups of rank~$2$, with the exception of the Ree groups
${}^2 F_4$. On the other hand, the theory developed in \cite{Mu99} allows to obtain twin root data by
integrating arbitrary Moufang foundations of $2$--spherical type. The groups obtained in this way are not
Kac-Moody groups whenever the foundation contains a Moufang octagon (which corresponds to a rank~$2$ parabolic
subgroup which is of type ${}^2 F_4$).
\end{itemize}

\vspace{2mm}

The conventions adopted throughout the rest of this section are the following: we let $(W,S)$ be a Coxeter
system with root system $\Phi$ and $(\Lambda, \{U_\alpha\}_{\alpha \in \Phi})$ be a twin root datum of type
$(W,S)$. We assume that each root group is finite and that $W$ is infinite. The associated twin buildings are
denoted $(\mathscr{B}_+, \mathscr{B}_-)$. The groups $\Aut(\mathscr{B}_\pm)$ are endowed with the compact-open
topology, which makes them locally compact totally disconnected second countable topological groups. We also
consider the subgroups $T, N, U_+, U_-$ of $\Lambda$ defined in \S\ref{sect:TRD}. The normal subgroup generated
by all root groups is denoted $\Lambda^\dagger$. If $g \in \Lambda$ fixes the building $\mathscr{B}_+$ it fixes
in particular the standard positive apartment and its unique opposite in $\mathscr{B}_-$ \cite[Lemma 2]{Ab95}
since $g$ preserves the twinning; therefore we have: $g \in T$ by \cite[Corollaire 3.5.4]{RemAst}. Moreover by
the Moufang parametrization (by means of root group actions) of chambers having a panel in the standard
apartment, such a $g$ must centralize each root group. This argument shows that the kernel of the action of
$\Lambda$ on $\mathscr{B}_+$ (resp. $\mathscr{B}_-$) is the centralizer $Z_\Lambda(\Lambda^\dagger)$ and is
contained in $T$.

\subsection{Topological completions: the building topology}\label{sect:BldgTopology}

For $\eps \in \{+, -\}$, let $\Lambda^{\mathrm{eff}}_\eps$ be the image of the natural homomorphism $\pi_\eps :
\Lambda \to \Aut(\mathscr{B}_\eps)$. Thus $\Lambda^{\mathrm{eff}}_+ \simeq \Lambda/Z_\Lambda(\Lambda^\dagger)
\simeq \Lambda^{\mathrm{eff}}_-$. The closure $\overline \Lambda^{\mathrm{eff}}_+ \leq \Aut(\mathscr{B}_+)$ is
the topological completion considered in \cite{RR06}. In the Kac-Moody case, another approach was proposed in
\cite{CarGar}, using the so-called weight topology; this allows to obtain completions of $\Lambda$ without
taking the effective quotient. However, the weight topology is defined using Kac-Moody algebras and, hence, does
not have an obvious substitute in the abstract framework considered here. Therefore, we propose the following.

\vspace{2mm}

For each non-negative integer $n$, let $K_\eps^n$ be the pointwise stabilizer in $U_\eps$ of the ball of
$\mathscr{B}_\eps$ centered at $c_\eps$ and of combinatorial radius $n$. Clearly $\bigcap_{n} K^n_{\eps} \subset
Z_\Lambda(\Lambda^\dagger) \subset T$ and, hence, $\bigcap_{n} K^n_{\eps} = \{1\}$ because $T \cap U_\eps =
\{1\}$ by \cite[Theorem~3.5.4]{RemAst}. Define a map $f_\eps : \Lambda \times \Lambda \to \R_+$ as follows:

\vspace{2mm}

\centerline{$f_\eps(g, h) = \left\{ \begin{array}{ll}
2  & \text{if } g\inv h \not \in U_\eps,\\
\exp(-\max \{ n \; | \; g\inv h \in K_\eps^n \}) &\text{if } g\inv h  \in U_\eps.
\end{array}
\right\}$.}

\vspace{2mm}

Since $K_\eps^n$ is a group for each $n$, it follows that $f_\eps$ is a left-invariant ultrametric distance on
$\Lambda$. We let $\overline \Lambda_\eps $ be the completion of $\Lambda$ with respect to this metric
\cite[II.3.7 Th\'eor\`eme 3 and III.3.4 Th\'eor\`eme 3.4]{BbkTG1-2-3-4};
this is the topological completion that we consider in this paper.

\vspace{2mm}

{\bf Definition.}~\it The so-obtained topology is called the {\rm building topology}~on $\overline\Lambda_\eps$.
\rm\vspace{2mm}

By left-invariance of the metric, replacing $U_\epsilon$ and $c_\epsilon$ by $\Lambda$--conjugates leads to the
same topology. Here is a summary of some of its basic properties; similar results hold with the signs $+$ and
$-$ interchanged.

\begin{prop}\label{prop:TopoCompletions}
We have the following:
\begin{itemize}
\item[(i)] The group $\overline \Lambda_+$ is locally compact and totally disconnected for the above topology.
It is second countable whenever $\Lambda$ is countable, i.e. whenever so is $T$.

\item[(ii)] The canonical map $\pi_+ : \Lambda \to \Lambda^{\mathrm{eff}}_+$ has a unique extension to a
continuous surjective open homomorphism $\overline \pi_+ :\overline \Lambda_+ \to \overline
\Lambda^{\mathrm{eff}}_+$.

\item[(iii)] The kernel of $\overline \pi_+$ is the discrete subgroup $Z_\Lambda(\Lambda^\dagger) < \overline
\Lambda_+$.

\item[(iv)] We have $\Stab_{\overline \Lambda_+}(c_+) \simeq  T \ltimes  \overline U_+$, where $\overline U_+$
denotes the closure of $U_+$ in $\overline \Lambda_+$.

\item[(v)] Every element $g \in \overline \Lambda_+$ may be written in a unique way as a product $g = u_+ n
u_-$, with $u_+ \in \overline U_+$, $n \in N$ and $u_- \in U_-$.

\item[(vi)] The sextuple $(\overline \Lambda_+, N, \overline U_+, U_-, T, S)$ is a refined Tits system, as
defined in \cite{KacPet}.
\end{itemize}
\end{prop}

\begin{remark}
It follows from (ii) and (iii) that the canonical map $\overline \Lambda_+ / Z_\Lambda(\Lambda^\dagger) \to
\overline \Lambda_+^{\mathrm{eff}}$ is an isomorphism of topological groups.
\end{remark}

\begin{proof} We start by noting that the restriction $\pi_+\!\mid_{U_+}$ is injective since $Z_\Lambda(\Lambda^\dagger) \cap U_\eps = \{1\}$ by \cite[Th\'eor\`eme 3.5.4]{RemAst}. Therefore, it follows from the definitions that $\pi_+ : U_+ \to \pi_+(U_+)$
is an isomorphism of topological groups.

\vspace{2mm}

We now prove (ii). Let $(\lambda_n)$ be a Cauchy sequence of elements of $\Lambda$. Let $n_0 \geq 0$ be such
that $f_+(\lambda_{n_0}, \lambda_n) \leq 1$ for all $n > n_0$. It follows that $\pi_+(\lambda\inv_{n_0}
\lambda_n)$ lies in the stabilizer of $c_+$ in $\Aut(\mathscr{B}_+)$, which is compact. This implies that
$\pi_+(\lambda_n)$ is a converging sequence in $\Aut(\mathscr{B}_+)$. Hence $\pi_+$ has a unique continuous
extension $\overline \pi_+ : \overline \Lambda_+ \to \overline \Lambda^{\mathrm{eff}}_+$ and it remains to prove
that $\overline \pi_+$ is surjective. By the preliminary remark above, it follows that $\overline \pi_+ :
\overline U_+ \to \overline{\pi_+(U_+)}$ is an isomorphism of topological groups. The surjectivity of $\overline
\pi_+$ follows easily since $\overline U_+$ is an open neighborhood of the identity. Finally, since $\overline
U_+$ contains a basis $\{\overline K_\eps^n\}$ of open neighborhoods of the identity, it follows that $\overline
\pi_+$ maps an open subset to an open subset.

\vspace{2mm}

(iv). The inclusion $T. \overline U_+ < \Stab_{\overline \Lambda_+}(c_+)$ is clear. Let $g \in
 \Stab_{\overline \Lambda_+}(c_+)$ and let $(\lambda_n)$ be a sequence in $\Lambda$
 such that
 $\displaystyle \lim_{n\to +\infty} \lambda_n = g$. Up to passing to a subsequence, we may -- and shall --
  assume that $\lambda_n \in
 \Stab_\Lambda(c_+)$ for all $n$. We know by \cite[\S~3.5.4]{RemAst} that $\Stab_\Lambda(c_+) \simeq
 T \ltimes U_+$. Hence each $\lambda_n$ has a unique writing $\lambda_n = t_n u_n$ with $t_n \in T$ and $u_n \in
 U_+$. Again, up to extracting a subsequence, we have $f_+(\lambda_1, \lambda_n) < 1$ for all $n$. In view of
 the semidirect decomposition $\Stab_\Lambda(c_+) \simeq
 T \ltimes U_+$, this implies that $t_n = t_1$ for all $n$. In particular, the sequence $(u_n)$ of elements of
 $U_+$ converges to $t_1\inv g$. This shows that $g \in  T.\overline U_+$ as desired. For every nontrivial
 element $t \in T$, we have $f_+(1, t) = 2$ because $T \cap U_+ = \{1\}$. On the other hand, for all $u \in
 \overline U_+$, we have $f_+(1, u) \leq 1$. Therefore, we have $T \cap \overline U_+ = \{1\}$.

\vspace{2mm}

(iii). The fact that $Z_\Lambda(\Lambda^\dagger)$ is discrete follows from $Z_\Lambda(\Lambda^\dagger) \cap U_+
= \{1\}$. Clearly we have $Z_\Lambda(\Lambda^\dagger) < \Ker (\overline \pi_+)$. We must prove the reverse
inclusion. Let $k \in \Ker (\overline \pi_+)$. By (iv), we have $k = tu$ for unique elements $t \in T$ and $u
\in \overline U_+$. Applying (iv) to the effective group $\overline \Lambda^{\mathrm{eff}}_+$, we obtain
$\pi_+(t) = 1$ and $\overline \pi_+(u) = 1$. Since the restriction of $\overline \pi_+$ to $\overline U_+$ is
injective by the proof of (ii), we deduce $u = 1$ and hence $k \in T < \Lambda$. Therefore $k \in
Z_\Lambda(\Lambda^\dagger)$ as desired.

\vspace{2mm}

(i). The building topology comes from a metric, so the Hausdorff group $\Lambda$ injects densely in its
completion $\overline \Lambda_+$  \cite[II.3.7 Corollaire]{BbkTG1-2-3-4} and the latter group is second
countable whenever $\Lambda$ is countable. It is locally compact because $\overline U_+$ is a compact open
subgroup by the proof of (ii). Furthermore, $\overline \pi_+$ annihilates the connected component of $\overline
\Lambda_+$ because $\overline \Lambda_+^{\mathrm{eff}}$ is totally disconnected. On the other hand, the kernel
of $\overline \pi_+$ is discrete by (iii). Hence $\overline \Lambda_+$ itself is totally disconnected.

\vspace{2mm}

(v). The group $U_-$ acts on $\mathscr{B}_+$ with the apartment $\mathscr{A}_+$ as a fundamental domain. The
group $N$ stabilizes $\mathscr{A}_+$ and acts transitively on its chambers. In view of (iv), it follows that
$\overline \Lambda_+ = \overline U_+ . N . U_-$. On the other hand, it follows easily from the definition of
$f_\eps$ that
$$\overline U_\eps = \{g \in \overline \Lambda_\eps \; | \; f_\eps(1,g) \leq 1\}.$$ Therefore, the uniqueness assertion
follows immediately from \cite[\S1.5.4]{RemAst} and the fact that $\Lambda \cap \overline U_+ = U_+$.

\vspace{2mm}

(vi). The main axiom of a refined Tits system is the property of assertion (v), which has just been proven. For
the other axioms to be checked, the arguments are the same as \cite[Proof of Theorem~1.C.(i)]{RR06}.\end{proof}

\subsection{Twin building lattices}\label{sect:TBLattice}
Let $q_{\min} = \min \{ |U_\alpha| : \alpha \in \Pi\}$, where $\Pi \subset \Phi$
is the set of simple roots. The following is an adaptation of the main result of \cite{RemCRAS}:

\begin{prop}\label{prop:lattices}
The image of the diagonal injection

\vspace{2mm}

\centerline{$\Lambda \to \overline \Lambda_+ \times \overline \Lambda_- : \lambda \mapsto (\lambda, \lambda)$}

\vspace{2mm}

is a discrete subgroup of $\overline \Lambda_+ \times \overline \Lambda_-$. It is an irreducible lattice if and
only if $W(1/q_{\min}) < + \infty$ and $Z_\Lambda(\Lambda^\dagger)$ is finite.
\end{prop}
\begin{proof}
Let $\Delta(\Lambda) = \{ (\lambda, \lambda) \; | \; \lambda \in \Lambda \} < \overline \Lambda_+ \times
\overline \Lambda_-$. The subgroup $\overline U_+ \times \overline U_- < \overline \Lambda_+ \times \overline
\Lambda_-$ is an open neighborhood of the identity. We have $\Delta(\Lambda) \cap (\overline U_+ \times
\overline U_-) = \Delta(U_+ \cap U_-)$. By \cite[\S3.5.4]{RemAst}, $ U_+ \cap U_- = \{1\}$. Thus
$\Delta(\Lambda)$ is discrete. The second assertion follows from the proofs of \cite[Proposition~5 and
Corollary~6]{Rem05}, which apply here without any modification: the only requirement is that $(\overline
\Lambda_+, N, \overline U_+, U_-, T, S)$ and $(\overline \Lambda_-, N, \overline U_-, U_+, T, S)$ be refined
Tits systems. This follows from Proposition~\ref{prop:TopoCompletions}(vi).
\end{proof}

Note that the group $Z_\Lambda(\Lambda^\dagger)$ may be arbitrarily large, since one may replace $\Lambda$ by
the direct product of $\Lambda$ with an arbitrary group; the root groups of $\Lambda$ also provide a twin root
datum for this direct product. However it is always possible to make the group $Z_\Lambda(\Lambda^\dagger)$
finite by taking appropriate quotients; note that if $\Lambda = \Lambda^\dagger$, then
$Z_\Lambda(\Lambda^\dagger) = Z(\Lambda)$ is abelian. Finally, since $Z_\Lambda(\Lambda^\dagger) < T$, it
follows that $Z_\Lambda(\Lambda^\dagger)$ is always finite when $\Lambda$ is a split or almost split Kac-Moody
group.

\vspace{2mm}

{\bf Definition.}~\it
If the twin root datum $(\Lambda, \{U_\alpha\}_{\alpha \in \Phi})$ is such that $W(1/q_{\min}) < +\infty$ and
$Z_\Lambda(\Lambda^\dagger)$ is finite, then $\Lambda$ is called a {\rm twin building lattice}.
\rm

\subsection{Structure of $\overline U_+$}

The first assertion of the following proposition was proved in \cite[Theorem~1.C(ii)]{RR06}.

\begin{prop}\label{prop:pronilpotency}
We have the following:
\begin{itemize}
\item[(i)] Assume that each root group $U_\alpha$ is a finite $p$-group. Then $\overline U_+$ is pro-$p$.

\item[(ii)] Assume that each root group $U_\alpha$ is solvable. Then $\overline U_+$ is pro-solvable.

\item[(iii)] Assume that each root group $U_\alpha$ is nilpotent. Then, for every prenilpotent set of roots
$\Psi \subset \Phi$, the group $U_\Psi = \la U_\alpha \; | \; \alpha \in \Psi \ra$ is nilpotent.
\end{itemize}
\end{prop}

\begin{remark}
One might expect that, under the assumption that all root groups are nilpotent, the group $\overline U_+$ is
pro-nilpotent. This is however not true in general. Counterexamples are provided by twin root data over ground
fields of mixed characteristics, constructed in \cite{RR06}. More precisely, consider a twin root datum
$(\Lambda, \{U_\alpha\}_{\alpha \in \Phi})$ of rank~$2$ with infinite Weyl group, such that the rank~one
subgroups are $\SL_2(\F_3)$ and $\SL_2(\F_4)$. The associated twin building is a biregular twin tree. Then the
$U_+$-action induced on the ball of combinatorial radius~$2$ centered at $c_+$ is not nilpotent: indeed, the
corresponding finite quotient of $U_+$ contains a subgroup isomorphic to the wreath product $(\Z/2\Z \times
\Z/2\Z) \wr \Z/3\Z$, which is not nilpotent.
\end{remark}

\begin{proof}
For (i), see \cite[1.C Lemma~1 p.198]{RR06}. The arguments given in [loc. cit.] can be immediately adapted to
provide a proof of (ii): the essential fact is that an extension of a solvable group (resp. $p$-group) by a
solvable group (resp. $p$-group) is again solvable.

\vspace{2mm}

(iii). A set of roots $\Psi$ is called \emph{nilpotent} if it is prenilpotent and if, moreover, for each pair
$\{\alpha, \beta \} \subset \Psi$ one has $[\alpha, \beta] \subset \Psi$. Since every prenilpotent set of roots
is contained in a nilpotent set (see \cite[\S1.4.1 and \S2.2.6]{RemAst}), it suffices to prove the assertion for
nilpotent sets. The proof is by induction on the cardinality of $\Psi$, the result being obvious when $\Psi$ is
a singleton. The elements of $\Psi$ can be ordered in a nibbling sequence $\alpha_1, \alpha_2, \dots, \alpha_n$;
hence the sets $\Psi_1 = \Psi \backslash \{\alpha_1\}$ and $\Psi_n = \Psi \backslash \{\alpha_n\}$ are nilpotent
[loc.~cit., \S1.4.1]. Furthermore, one has $[U_{\alpha_1}, U_{\Psi_1}] \leq U_{\Psi_1}$ and $[U_{\alpha_n},
U_{\Psi_n}] \leq   U_{\Psi_n}$ as a consequence of (TRD1). Therefore, the subgroups $U_{\Psi_1}$ and
$U_{\Psi_n}$ are normal in $U_\Psi$, and are nilpotent by the induction hypothesis. It follows that $U_\Psi$ is
nilpotent \cite[Theorem~10.3.2]{Hall}. This part of the proof does not require that the root groups be finite.
\end{proof}

\section{Further properties of topological completions}\label{s - AptTopology}

In this section, we introduce a property of fixed points of root subgroups of a group $\Lambda$ endowed with a
twin root datum; this property is called (FPRS). We first provide sufficient conditions which ensure that this
property holds for any split or almost split Kac-Moody group, as well as for all exotic twin building lattices
mentioned in Sect.~\ref{ss - exoticTRD}. We then show that (FPRS) implies that the topological completion
$\overline{\Lambda}_+$ is topologically simple (modulo the kernel of the action on the building, see
Proposition~\ref{prop:TopoSimplicity}). Property~(FPRS) will be used again below, as a sufficient condition
implying that any nontrivial continuous homomorphism whose domain is $\overline \Lambda_+$ is proper
(Theorem~\ref{th:proper}).

\vspace{2mm}

Throughout this section, we let $\Lambda$ be a group endowed with a twin root datum $\{U_\alpha\}_{\alpha \in
\Phi}$ of type $(W,S)$ and let $(\mathscr{B}_+, \mathscr{B}_-)$ be the associated twin buildings.

\subsection{Fixed points of root subgroups}\label{sect:FPRS}

For any subgroup $\Gamma \leq G$, we define $r(\Gamma)$ to be the  supremum of the set of all non-negative real
numbers $r$ such that  $\Gamma$ fixes pointwise the (combinatorial) closed ball $B(c_+, r)$ of (combinatorial)
radius $r$ centered at $c_+$. In the present subsection, we consider the following condition:

\begin{description}
\item[\textbf{(FPRS)}] \emph{Given any  sequence of roots $(\alpha_n)_{n \geq 0}$  of $\Phi(\mathscr{A})$ such
that $\displaystyle \lim_{n \to +\infty} d(c_+, \alpha_n) = +\infty$, we have: $\displaystyle \lim_{n \to
+\infty} r(U_{-\alpha_n}) = +\infty$.}
\end{description}

\begin{remark}
This property can be seen as a non-quantitative generalization of \cite[7.4.33]{BruhatTits72}.
\end{remark}

In other words, this means that if the sequence of roots $(\alpha_n)_{n \geq 0}$ is such that $\displaystyle
\lim_{n \to +\infty} d(c_+, \alpha_n) = +\infty$, then the sequence of root subgroups $(U_{-\alpha_n})_{n \geq
0}$ tends to the identity in the building topology. The purpose of this subsection is to establish sufficient
conditions on the twin root datum $(\Lambda, (U_\alpha)_{\alpha \in \Phi})$ which ensure that property (FPRS)
holds. To this end, we will need the following conditions:

\begin{description}
\item[\textbf{(PP)}] \emph{For any prenilpotent pair of roots $\{\alpha, \beta\}$ such that $\la r_\alpha,
r_\beta \ra$ is infinite, either $[U_\alpha, U_\beta] = \{1\}$ or there exists a root $\phi$ such that
$r_\phi(\alpha) = -\beta$,
$[U_\alpha, U_\beta] \leq U_\phi$ and $[U_\alpha, U_\phi] = \{1\} = [U_\beta, U_\phi]$.}\\

\item[\textbf{($2$-sph)}] \emph{The Coxeter system $(W,S)$ is  $2$-spherical
and $G$ possesses no critical rank $2$  subgroup.}\\
This means that any pair of elements of $S$ generates a finite group  and moreover that for any pair $\{\alpha,
\beta\} \subset \Pi$, the group $X_{\alpha, \beta}$ generated by the  four root groups $U_{\pm\alpha}$,
$U_{\pm\beta}$, divided by its  center, is not isomorphic to any of the groups $B_2(2), G_2(2), G_2(3)$  or
${}^2 F_4(2)$.
\end{description}

\vspace{2mm}
The main result of this section is the following:

\begin{prop}\label{prop:FPRS}
Assume that the twin root datum $(\Lambda, (U_\alpha)_{\alpha \in \Phi})$ satisfies (PP) or ($2$-sph) or that
$\Lambda$ is a split or almost split Kac-Moody group. Then the property (FPRS) holds.
\end{prop}

\begin{remark}
The exotic examples of twin root data mentioned in Sect.~\ref{ss - exoticTRD} (I) and (II) satisfy condition
(PP). In fact, they satisfy even the stronger condition that all commutation relations are trivial: for any
prenilpotent pair of distinct roots $\{\alpha, \beta\}$, one has $[U_\alpha, U_\beta] = \{1\}$. Moreover, the
examples of type (III) satisfy ($2$-sph). However, it was communicated to us by B.~M\"uhlherr that there exists
an example of a group endowed with a twin root datum, which does not satisfy condition (FPRS). In this example,
whose construction is nontrivial, the Weyl group is the free Coxeter group of rank~$3$ (i.e. a free product of
$3$~copies of the group of order~$2$) and the ground field is $\F_2$.
\end{remark}

The proof of Proposition~\ref{prop:FPRS} splits into a sequence of lemmas which we prove separately.

\begin{lem}\label{lem:FPRS}
Assume that (PP) holds. For each integer $n \geq 0$, each root $\alpha \in \Phi(\mathscr{A})$ and each chamber
$c \in \mathscr{A}_+$, if $d(c, \alpha) \geq \frac{4^{n+1}-1}{3}$, then $U_{-\alpha}$ fixes $B(c,n)$ pointwise.
In particular (FPRS) holds.
\end{lem}
\begin{proof}
We work by induction on $n$. If $d(c, \alpha) \geq 1$, then $c \not \in \alpha$ whence $c \in -\alpha$. In
particular $c$ is fixed by $U_{-\alpha}$. Thus the desired property holds for $n=0$.

\vspace{2mm}

Assume now that $n>0$ and let $\alpha$ be a root such that $d(c, \alpha) \geq \frac{4^{n+1}-1}{3}$. By
induction, the group $U_{-\alpha}$ fixes the ball $B(c, n-1)$ pointwise. Furthermore, if $c'$ is a chamber
contained in $\mathscr{A}_+$ and adjacent to $c_+$, then $d(c', \alpha) \geq d(c,\alpha)-1$; therefore, the
induction hypothesis also implies that $U_{-\alpha}$ fixes the ball $B(c', n-1)$ pointwise.

\vspace{2mm}

Let now $x$ be a chamber at distance $n$ from $c$. Let $c_0 = c, c_1, \dots, c_n= x$ be a minimal gallery
joining $c$ to $x$. We must prove that $U_{-\alpha}$ fixes $x$. If $c_1$ is contained in $\mathscr{A}_+$ then we
are done by the above. Thus we may assume that $c_1$ is not in $\mathscr{A}_+$. Let $c'$ be the unique chamber
of $\mathscr{A}_+$ such that $c, c_1$ and $c'$ share a panel. Let $\beta \in \Phi(\mathscr{A})$ be one of the
two roots such that
the wall $\partial \beta$ separates $c$ from $c'$. Up to replacing $\beta$ by its opposite if necessary, we may
-- and shall -- assume by \cite[Proposition~9]{Ti88/89} that the pair $\{-\alpha, \beta\}$ is prenilpotent. Let
$u \in U_\beta$ be the (unique) element such that $u(c_1)$ belongs to $\mathscr{A}_+$; thus we have $u(c_1)= c$
or $c'$. Since $u(c_1), u(c_2), \dots, u(c_n)$ is a minimal gallery, it follows that $u(x)$ is contained in
$B(c, n-1) \cup B(c', n-1)$.

\vspace{2mm}

There are three cases.

\vspace{2mm}

Suppose first that $[U_{-\alpha}, U_\beta]=\{1\}$. For any $g \in U_{-\alpha}$, we have $g = u\inv g u$ whence
$g(x) = u\inv g u(x) = x$ because $g \in U_{-\alpha}$ fixes $B(c, n-1) \cup B(c', n-1)$ pointwise.

\vspace{2mm}

Suppose now that $[U_{-\alpha}, U_\beta] \neq \{1\}$ and that $\la r_\alpha, r_\beta \ra$ is infinite. By
property (PP) there exists a root $\phi \in \Phi(\mathscr{A})$ such that $[U_{-\alpha}, U_\beta] \leq U_\phi$
and $r_\phi(\alpha) = \beta$. Let $y_0=c, y_1, \dots, y_k$ be a gallery of minimal possible length joining $c$
to a chamber of $-\phi$. Thus we have $y_k \in -\phi$, $y_{k-1} \in \phi$ and $k = d(c, -\phi)$. Since
$r_\phi(\beta)=\alpha$ and since either $c$ or $c'$ belongs to $\beta$, it follows from considering the gallery
$$c=y_0, \dots, y_{k-1}, y_k =r_\phi(y_{k-1}), r_\phi(y_{k-2}), \dots, r_\phi(c), r_\phi(c')$$
of length $2k$, that $d(c, \alpha) \leq 2k$, whence $d(c, -\phi) \geq \frac{1}{2}d(c, \alpha) \geq
\frac{4^{n+1}-1}{6} > \frac{4^n-1}{3}$. A similar argument shows that $d(c', -\phi) \geq \frac{1}{2}d(c',
\alpha) \geq \frac{4^{n+1}-4}{6}> \frac{4^n-1}{3}$. Therefore, the induction hypothesis shows that $U_\phi$
fixes $B(c, n-1) \cup B(c', n-1)$ pointwise. Now, for any $g \in U_{-\alpha}$, we have $g(x) = [g, u\inv] u\inv
g u (x) = [g, u\inv](x)$ because $g \in U_{-\alpha}$ fixes $B(c, n-1) \cup B(c', n-1)$ pointwise. By (PP), the
commutator $[g, u\inv]$ commutes with $u$ and, hence, we have $[g, u\inv](x) = u\inv [g, u\inv]u(x) = x$ because
$[g, u\inv] \in U_\phi$ fixes $B(c, n-1) \cup B(c', n-1)$ pointwise.

\vspace{2mm}

Suppose finally that $[U_{-\alpha}, U_\beta] \neq \{1\}$ and that $\la r_\alpha, r_\beta \ra$ is finite. This
implies that the pairs $\{-\alpha, \beta\}$ and $\{-\alpha, -\beta\}$ are both prenilpotent. Therefore, up to
replacing $\beta$ by its opposite if necessary, we may -- and shall -- assume that $c \not \in \beta$, whence
$u(c_1) = c$. Note that $\la r_\alpha, r_\beta \ra$ is contained in a rank~$2$ parabolic subgroup $P$ of $W$.
Since any such subgroup is the Weyl group of a Levi factor of rank~$2$ of $\Lambda$, which is itself endowed
with a twin root datum of rank~$2$, it follows from \cite[Theorem~17.1]{TW02} that $P$ is of order at most~$16$.
Let $[-\alpha, \beta]= \{\gamma \in \Phi(\mathscr{A}) \; | \; (-\alpha) \cap \beta \subseteq \gamma, \ \alpha
\cap (-\beta) \subseteq -\gamma \}$; thus $[-\alpha, \beta]$ has at most $8$~elements because for every $\gamma
\in [-\alpha, \beta]$, the reflection $r_\gamma$ belongs to $P$. Order the elements of $[\beta, -\alpha]$ in a
natural cyclic order: $[\beta, -\alpha] =\{\beta = \beta_0, \beta_1, \dots, \beta_m = -\alpha\}$; this means
that $r_{\beta_i}(\beta_{i-1}) = \beta_{i+1}$ for $i = 1, \dots, k-1$. Such an ordering does exist because the
group $\la  r_{\gamma} \; | \; \gamma \in [\beta, -\alpha] \ra$ is (finite) dihedral. Let $c=y_0, y_1, \dots,
y_k$ be a gallery of minimal possible length joining $c$ to a chamber of $-\beta_1$. Thus we have $y_k \in
-\beta_1$, $y_{k-1} \in \beta_1$ and $k = d(c, -\beta_1)$. Since $r_{\beta_1}(\beta)=-\beta_2$ and since $c'$
belongs to $\beta$, it follows from considering the gallery
$$c=y_0, \dots, y_{k-1}, y_k =r_{\beta_1}(y_{k-1}), r_{\beta_1}(y_{k-2}), \dots, r_{\beta_1}(c), r_{\beta_1}(c')$$
of length $2k$, that $d(c, -\beta_2) \leq 2k$. A straightforward induction yields $d(c, -\beta_i) \leq ik$ for
$i = 1, \dots, m$. In particular, we have $d(c, \alpha) \leq m k = m. d(c, -\beta_1)$. Recall that  $m +1$ is
the cardinality of $[-\alpha, \beta]$.

\vspace{2mm}

We may now choose a natural cyclic order
$[-\beta, -\alpha] = \{-\beta = \beta'_0, \beta'_1, \dots, \beta'_{m'} = -\alpha\}$ and
repeat the same arguments with $c$ replaced by $c'$, $\beta$ replaced by $-\beta$ and each $\beta_i$
replaced by $r_\beta(\beta_i)$. This yields $d(c', -\beta'_i) \leq i d(c', -\beta'_1)$ for each $i$.
Note that $d(c, -\beta_1) = d(r_\beta(c), -r_\beta(\beta_1)) = d(c', -\beta'_1)$. We
obtain that $d(c, \alpha)-1 \leq d(c', \alpha) \leq m'. d(c, -\beta_1)$,
where $m' +1$ is the cardinality of
$[-\alpha, -\beta]$. Observe now that $m + m'= \frac{|P|}{2} \leq 8$. In particular, we have $\min\{m, m'\} \leq
4$. Therefore, we deduce from the inequalities above that for each $i = 1, \dots, m$, we have
$$d(c, -\beta_i) \geq d(c, -\beta_1) \geq \frac{d(c, \alpha)-1}{4} \geq \frac{4^n-1}{3}.$$ By the induction
hypothesis, it follows that for each $\gamma \in ]{-\alpha}, \beta[ = [-\alpha, \beta] \backslash \{ -\alpha,
\beta\} $, the root subgroup $U_{\gamma}$ fixes the ball $B(c, n-1)$ pointwise.

\vspace{2mm}

Now, for any $g \in U_{-\alpha}$, we have $g(x) = [g, u\inv] u\inv g u (x) = [g, u\inv](x)$ because $g \in
U_{-\alpha}$ fixes $B(c, n-1)$ pointwise. Moreover, we have $[g, u\inv] \in U_{]{-\alpha}, \beta[} = \la
U_\gamma\; | \; \gamma \in ]{-\alpha}, \beta[ \ra$ by (TRD1). Therefore, the commutator $[g, u\inv]$ fixes
$u(x)$ and we have
$$g(x) =  [g, u\inv](x) = \big([[g, u\inv], u\inv]\big) u\inv[g, u\inv] u(x) = [[g, u\inv], u\inv](x).$$
Repeating the argument $m$ times successively, we obtain $g(x) =[\dots [[g, u\inv], u\inv], \dots, u\inv ](x)$
where the commutator is iterated $m$ times. By (TRD1), we have
$[\dots [[g, u\inv], u\inv], \dots, u\inv ] \in U_{]\beta_1, \beta[}$,
which is trivial since $]\beta_1, \beta[$ is empty. Therefore, we deduce finally that $g$ fixes $x$, as desired.
\end{proof}

\begin{lem}\label{lem:splitKM=>PP}
Suppose that $\Lambda$ is a split Kac-Moody groups and that $(U_\alpha)_{\alpha \in \Phi}$ is its natural system
of root subgroups. Then the twin root datum $(\Lambda, (U_\alpha)_{\alpha \in \Phi})$ satisfies (PP).
\end{lem}
\begin{proof}
This follows by combining \cite[Theorem~2]{Morita87} with some results from \cite{BilligPia95}
(see also \cite[Sect.~3.2]{TitsJA}).
In order to be more precise, we freely use the notation and terminology of these references in the present proof. In
particular, we use the `linear' root system of the Lie algebra associated with the Kac-Moody group $\Lambda$,
instead of the `abstract' root system introduced above and which is appropriate to the case of general twin root
data. A comprehensive introduction to linear root systems can be found in \cite[Chapter~5]{MP95}.

\vspace{2mm}

Commutation relations in split Kac-Moody groups are described precisely by \cite[Theorem~2]{Morita87}. Combining
the latter result together with \cite[Proposition~1]{BilligPia95}, one sees easily that if $\{\alpha, \beta\}$
is a prenilpotent pair such that $\la r_\alpha, r_\beta\ra$ is infinite and $[U_\alpha, U_{\beta}] \neq \{1\}$,
then $[U_\alpha, U_{\beta}] \leq U_{\alpha +\beta}$ and  the $\alpha$-string through $\beta$ is of length~$\geq
5$ and contains exactly $4$~real roots, which are $\beta - \la \beta, \alpha^{\vee} \ra\alpha$, $\beta - (\la
\beta, \alpha^{\vee} \ra - 1)\alpha$, $\beta$ and $\beta + \alpha$. In particular, $\beta + 2\alpha$ is not a
root, whence $\{-\alpha, \alpha + \beta\}$ is $W$-conjugate to a Morita pair by
\cite[Proposition~3(i)]{BilligPia95}. In particular, we have $\la -\alpha, (\alpha + \beta)^\vee \ra = -1$ by
\cite[Proposition~2]{BilligPia95}. We deduce that $r_{\alpha + \beta}(-\alpha)=-\alpha + \la -\alpha, (\alpha +
\beta)^\vee \ra(\alpha +\beta )= \beta$. Finally, since $2\alpha + \beta$ is not a root, it follows from
\cite[Theorem~2]{Morita87} that $[U_\alpha, U_{\alpha + \beta}] = \{1\}$. Hence property (PP) holds, as desired.
\end{proof}

\begin{lem}\label{lem:AlmostSplitKM=>FPRS}
Suppose that $\Lambda$ is an almost split Kac-Moody groups and that $(U_\alpha)_{\alpha \in \Phi}$ is its
natural system of root subgroups. Then the twin root datum $(\Lambda, (U_\alpha)_{\alpha \in \Phi})$ satisfies
(FPRS).
\end{lem}
\begin{proof}
Our reference on almost split Kac-Moody groups is \cite[Chapitres~11--13]{RemAst}. Let $\K$ be the ground field
of $\Lambda$, let $\K_{\mathrm{s}}$ be a separable closure of $\K$, let $\Gamma =
\mathrm{Gal}(\K_{\mathrm{s}}/\K)$ and let $\widetilde{\Lambda}$ be a split Kac-Moody group over
$\K_{\mathrm{s}}$, such that $\Lambda$ is the fixed point set of a $\Gamma$-action on $\widetilde{\Lambda}$. We
henceforth view $\Lambda$ as a subgroup of $\widetilde{\Lambda}$. We denote by $(\widetilde{U}_\alpha)_{\alpha
\in \widetilde{\Phi}}$ the natural system of root subgroups of $\widetilde{\Lambda}$ and by
$(\widetilde{\mathscr{B}}_+, \widetilde{\mathscr{B}}_-)$ the twin building associated with the twin root datum
$(\widetilde{\Lambda}, (\widetilde{U}_\alpha)_{\alpha \in \widetilde{\Phi}})$. By
\cite[Th{\'e}or{\`e}me~12.4.4]{RemAst}, the twin building $(\mathscr{B}_+, \mathscr{B}_-)$ is embedded in a
$\Lambda$-equivariant way in $(\widetilde{\mathscr{B}}_+, \widetilde{\mathscr{B}}_-)$, as the fixed point set of
$\Gamma$-action on $(\widetilde{\mathscr{B}}_+, \widetilde{\mathscr{B}}_-)$. This embedding maps chambers of
$(\mathscr{B}_+, \mathscr{B}_-)$ to $\K$-chambers of $(\widetilde{\mathscr{B}}_+, \widetilde{\mathscr{B}}_-)$,
which are minimal $\Gamma$-invariant spherical residues. Let $r$ be the rank of such a spherical residue. Two
$\K$-chambers are adjacent (as chambers of $(\mathscr{B}_+, \mathscr{B}_-)$) if they are contained in a common
spherical residue of rank~$r+1$ and either coincide or are opposite in that residue. This shows that bounded
subsets of $(\mathscr{B}_+, \mathscr{B}_-)$ are also bounded in $(\widetilde{\mathscr{B}}_+,
\widetilde{\mathscr{B}}_-)$ and, moreover, that every ball of large radius in $(\widetilde{\mathscr{B}}_+,
\widetilde{\mathscr{B}}_-)$ which is centered at a point of a $\K$-chamber contains a ball of large radius of
$(\mathscr{B}_+, \mathscr{B}_-)$.

\vspace{2mm}

Let $(\alpha_n)_{n \geq 0}$ be a sequence of roots of $\Phi=\Phi(\mathscr{A})$ such that $d(c_+, \alpha_n)$
tends to infinity with $n$. We must prove that $r(U_{-\alpha_n})$ also tends to infinity with $n$. We choose the
base chamber $\widetilde{c}_+$ of $(\widetilde{\mathscr{B}}_+, \widetilde{\mathscr{B}}_-)$ such that it is
contained in the $\K$-chamber $c_+$, and denote by $\widetilde{r}(H)$ the supremum of the radius of a ball
centered at $\widetilde{c}_+$ which is pointwise fixed by $H$. In view of the preceding paragraph, it suffices
to show that $\widetilde{r}(U_{-\alpha_n})$ tends to infinity with $n$. To this end, we will use the fact that
$(\widetilde{\Lambda}, (\widetilde{U}_\alpha)_{\alpha \in \widetilde{\Phi}})$ satisfies property (FPRS) by
Lemmas~\ref{lem:FPRS} and~\ref{lem:splitKM=>PP}. Let $\beta \in \Phi$ be a $\K$-root and consider the root
subgroup $U_\beta$. Let $x, y$ be two adjacent $\K$-chambers such that $\beta$ contains $x$ but not $y$. Let
$\widetilde{\Phi}(\beta)$ be the set of ($\K_{\mathrm{s}}$-)roots containing $x$ but not $y$; it is independent
of the choice of $x$ and $y$. Furthermore $\widetilde{\Phi}(\beta)$ is a prenilpotent subset of
$\widetilde{\Phi}$ and $U_\beta \subset \widetilde{U}_{\widetilde{\Phi}(\beta)} = \la \widetilde{U}_\gamma \; |
\; \gamma \in \widetilde{\Phi}(\beta) \ra$ by \cite[\S12.4.3]{RemAst}. Therefore, in order to finish the proof,
it suffices to show that $\min \{ d(\widetilde{c}_+, \gamma) \; | \; \gamma \in \widetilde{\Phi}(\alpha_n)\}$
tends to infinity with $n$.

\vspace{2mm}

Assume for a contradiction that this is not the case. Then there exists a subsequence $(\alpha_{n_j})_{ j \geq
0}$ and an element $\gamma_j \in \widetilde{\Phi}(\alpha_{n_j})$ such that $d(\widetilde{c}_+, \gamma_j)$ is a
bounded function of $j$. Since the apartments are locally finite, it follows that, up to extracting a
subsequence, we may -- and shall -- assume that $\gamma_j$ is constant. On the other hand, we claim (and prove
below) that if $\alpha, \alpha'$ are two distinct $\K$-roots, then the sets $\widetilde{\Phi}(\alpha)$ and
$\widetilde{\Phi}(\alpha')$ are disjoint; this implies that the sequence $(\alpha_{n_j})_{ j \geq 0}$ is
constant, contradicting the fact that $d(c_+, \alpha_{n_j})$ tends to infinity with $j$. It remains to prove the
claim. This is most easily done using the notion of (combinatorial) projections in buildings. Let $\pi$ (resp.
$\pi'$) be a $\K$-panel stabilized by the $\K$-reflection $r_\alpha$ (resp. $r_{\alpha'}$). Thus $\pi$ and
$\pi'$ are spherical residues of $\widetilde{\mathscr{B}}_+$ of rank $r+1$. Note that $\proj_\pi(\pi')$ (resp.
$\proj_{\pi'}(\pi)$) contains a $\K$-chamber. Furthermore, given any $\gamma \in \widetilde{\Phi}(\alpha) \cap
\widetilde{\Phi}(\alpha')$, the ($\K_{\mathrm{s}}$-)reflection $r_\gamma$ stabilizes both $\pi$ and $\pi'$, but
it does not stabilize any $\K$-chamber. Therefore $\proj_\pi(\pi')$ (resp. $\proj_{\pi'}(\pi)$) cannot be
reduced to a single $\K$-chamber; since $\proj_\pi(\pi')$ (resp. $\proj_{\pi'}(\pi)$) is a sub-residue of $\pi$
(resp. $\pi'$), it must be of rank $r+1$, which yields $\proj_\pi(\pi')=\pi$ (resp. $\proj_{\pi'}(\pi)=\pi'$).
In other words the $\K$-panels $\pi$ and $\pi'$ are \emph{parallel}. Since $r_\alpha$ stabilizes any $\K$-panel
which is parallel to $\pi$ by \cite[Proposition~2.7]{CM04}, this implies that the $\K$-reflections $r_\alpha $
and $r_{\alpha'}$ stabilize a common $\K$-panel. Since $\alpha \neq \alpha'$, we deduce that $\alpha =
-\alpha'$, which implies that $\widetilde{\Phi}(\alpha) = -\widetilde{\Phi}(\alpha') = \{-\gamma \; | \; \gamma
\in \widetilde{\Phi}(\alpha')\}$. In particular, the sets $\widetilde{\Phi}(\alpha)$ and
$\widetilde{\Phi}(\alpha')$ are disjoint, a contradiction.
\end{proof}

\begin{lem}\label{lem:2sph}
Suppose that $\Lambda$ satisfies condition {\rm ($2$-sph)} and that each root group is finite. Then property
(FPRS) holds.
\end{lem}
\begin{proof}
As before, we denote by $(W,S)$ be the Coxeter system consisting of the Weyl group $W$ together with its
canonical generating set $S$. If $(W,S)$ is not of irreducible type, then the buildings $\mathscr{B}_+$ and
$\mathscr{B}_-$ split into direct products of irreducible components, and it is easy to see that checking (FPRS)
for the $\Lambda$-action on $\mathscr{B}_+$ is equivalent to checking (FPRS) for the induced action on each
irreducible component. We henceforth assume that $(W,S)$ is of irreducible type. If $W$ is finite, then there is
no sequence of roots $(\alpha_n)_{n \geq 0}$ such that $d(c_+, \alpha_n)$ tends to infinity with $n$. Assume now
that $W$ is infinite; in particular $(W,S)$ is of rank~$\geq 3$. Consider a sequence of roots $(\alpha_n)_{n
\geq 0}$ such that $d(c_+, \alpha_n)$ tends to infinity with $n$. We must prove that $r(U_{-\alpha_n})$ tends to
infinity with $n$.

\vspace{2mm}

Given any two basis roots $\alpha, \alpha'$, there exists a sequence of basis roots $\alpha = \alpha_0,
\alpha_1, \dots, \alpha_k = \alpha'$ such that $r_{\alpha_{i-1}}$ does not commute with $r_{\alpha_i}$ for $i =
1, \dots, k$ because $(W,S)$ is irreducible. This implies that each rank two subgroup $X_{\alpha_{i-1},
\alpha_i} = \la U_{\pm \alpha_{i-1}}, U_{\pm \alpha_{i}} \ra$ is endowed with a twin root datum of irreducible
spherical type and rank~$2$, since $(W,S)$ is $2$-spherical. By the classification of such groups \cite{TW02},
it follows that $U_{\alpha_{i-1}}$ and $U_{\alpha_{i}}$ are finite $p_i$-groups for some prime $p_i$. Since this
is true for each $i$, we have $p_i = p_{i-1}$, whence the sequence $(p_i)$ is constant. This shows that there
exists a prime $p$ such that each root group is a finite $p$-group. By Proposition~\ref{prop:pronilpotency}(i),
it follows that $\overline{U}_+$ is pro-$p$, whence pro-nilpotent. In particular, the descending central series
$(U_+^{(n)})_{n \geq 0}$ tends to the identity in the building topology when $n \to +\infty$.
Therefore, in order to finish the proof, it suffices to show that $\displaystyle \lim_{n \to +\infty}
q(-\alpha_n) = +\infty$, where $q(\alpha) =  \max \{ n \geq 0 : U_\alpha \leq U_+^{(n)} \}$ for each positive
root $\alpha \in \Phi_+$.

\vspace{2mm}

For each integer $n \geq 0$, we set $\Phi_+^n = \{\alpha \in \Phi_+ : q(\alpha) = n\}$. We claim that each
$\Phi_+^n$ is finite. By assumption, for every $n$ there exists $n'$ such that the set $\{\alpha_j : j \geq
n'\}$ does not contain any element of $\Phi_+^n$, the desired assertion follows from the claim. In order to
prove the claim, we proceed by induction on $n$. We first need to recall a consequence of condition ($2$-sph).

\vspace{2mm}

A pair $\{\alpha, \beta\} \subset \Phi_+$ is called  \textit{fundamental} if the following conditions hold:
\begin{description}
\item[\textbf{(FP1)}] The group  $\la r_\alpha, r_\beta \ra$ is finite.

\item[\textbf{(FP2)}] For each $\gamma \in \Phi_+$ such that the group  $\la r_\alpha, r_\beta, r_\gamma \ra$ is
dihedral, we have $\gamma \in [\alpha, \beta]$. In other words, this  means that the pair $\{\alpha, \beta\}$ is
a basis of the root subsystem it generates.
\end{description}
We have the following:
\begin{itemize}
\item[(a)] Let $\{\alpha, \beta\} \subset \Phi_+$ be a fundamental  pair. Then, for all $\gamma \in ]\alpha,
\beta[$, we have $U_\gamma \leq [U_\alpha, U_\beta]$ by \cite[Proposition~7]{Ab95}.

\item[(b)] Let $\gamma \in \Phi_+$ be a root such that  $d(c_+,-\gamma)>1$. Then there exists a fundamental pair
$\{\alpha, \beta\} \subset \Phi_+$ such that $\gamma \in ]\alpha,  \beta[$. This follows from
\cite[Lemma~1.7]{BH94} together with  the fact that $(W,S)$ is $2$-spherical.
\end{itemize}

We now prove by induction on $n$ that $\Phi_+^n$ is finite. The set $\Phi_+^0$ coincides with $\Pi$. Indeed, for
each simple root $\alpha \in \Pi$, the group $U_\alpha$ fixes $c_+$ but acts non-trivially on the chambers
adjacent to $c_+$. Since on the other hand, the derived group $U_+^{(1)}$ fixes the ball  $B(c_+,1)$ pointwise,
we deduce that $U_\alpha$ is not contained in $U_+^{(1)}$, whence $q(\alpha)=0$. Thus $\Pi \subset \Phi_+^0$.
Conversely, if $\alpha \in \Phi_+$ does not belong to $\Pi$, then  $d(c_+, -\alpha)>1$ and property~(a) implies
that $U_\alpha \leq U_+^{(1)}$. Thus $q(\alpha) \geq 1$ and $\alpha \not \in \Phi_+^0$. This shows that
$\Phi_+^0 = \Pi$. In particular $\Phi_+^0$ is finite.

\vspace{2mm}

Let now $n\geq 1$ and assume that $\Phi_+^k$ is finite for all $k < n$. We must prove that $\Phi_+^n$ is finite.
Let us enumerate its elements: $\Phi_+^n = \{\gamma_1, \gamma_2, \dots  \}$. Since $n\geq 1$ and since $\Phi_+^0
=\Pi$, we have $d(c_+,  -\gamma_i)>1$ for all $i \geq 1$. Hence, by property~(b), for each $i$ there is  a
fundamental pair $\{\alpha_i, \beta_i\}$ such that $\gamma_i \!\in\,  ]\alpha_i, \beta_i[$. By property~(b),
this implies $U_{\gamma_i} <  [U_{\alpha_i}, U_{\beta_i}]$. Therefore, we have $n = q(\gamma) >
\max\{q(\alpha_i), q(\beta_i)\}$. In particular:

\vspace{2mm}

\centerline{$\displaystyle\bigcup_{i>0}\{\alpha_i, \beta_i\} \subset  \bigcup_{k =0}^{n-1} \Phi_+^k$.}

\vspace{2mm}

The set $\displaystyle\bigcup_{i>0} \ \{\alpha_i, \beta_i\}$ is thus  finite. By the definition of the
$\gamma_i$'s, we have \vspace{2mm}

\centerline{$\displaystyle\Phi_+^n \subset \bigcup_{i>0} \ ]\alpha_i,  \beta_i[$.}

\vspace{2mm}

Since each interval $]\alpha_i, \beta_i[$ is finite  \cite[2.2.6]{RemAst}, this shows that $\Phi_+^n$ is finite.
\end{proof}

\subsection{Density of the commutator subgroup}

As before, $(\Lambda, \{U_\alpha\}_{\alpha \in \Phi})$ is a twin root datum of type $(W,S)$ and $(\mathscr{B}_+,
\mathscr{B}_-)$ is the associated twin buildings. We assume moreover, for the rest of this section, that all
root groups are finite.

\begin{lem}\label{lem:density}
Assume that property (FPRS) holds, that the Weyl group $W$ is infinite and that the associated Coxeter system
$(W,S)$ is irreducible. If $\Lambda$ is generated by its root subgroups, then the commutator subgroup
$[\overline \Lambda_+, \overline \Lambda_+]$ is dense in $\overline \Lambda_+$.
\end{lem}

\begin{remark}
If each rank one subgroup $X_\alpha = \la U_\alpha \cup U_{-\alpha} \ra$ of $\Lambda$ is perfect and if
$\Lambda$ is generated by its root subgroups, then $[\overline \Lambda_+, \overline \Lambda_+] \supset \Lambda$
and, hence, $[\overline \Lambda_+, \overline \Lambda_+]$ is dense in $\overline \Lambda_+$. However, there are
many examples of groups endowed with a twin root datum satisfying (FPRS) but whose rank one subgroups are not
perfect, e.g. Kac-Moody groups over $\F_2$ or $\F_3$, or twin building lattices as in Sect.~\ref{ss -
exoticTRD}(II) where the rank one subgroups may be solvable.
\end{remark}

\begin{proof}
Let $\varphi : \overline \Lambda_+ \to G$ be a continuous homomorphism to an abelian topological group $G$. Let
$\Pi$ be the standard root basis of $\Phi$, where $\Phi$ is the root system of $(W,S)$ indexing the twin root
datum of $\Lambda$. For each $\alpha \in \Pi$, let $X_\alpha = \la U_\alpha \cup U_{-\alpha} \ra$. Applying
Lemma~\ref{lem:Rank1} below  to $\varphi|_{X_\alpha}$, it follows  that $\varphi (U_\alpha )=
\varphi(U_{-\alpha})$ for all $\alpha \in \Pi$.

\vspace{2mm}

Assume by contradiction that $\varphi$ is nontrivial. Since $\Lambda = \la U_\alpha \; | \; \alpha \in \Phi \ra
= \la U_\alpha \; | \; \alpha \in \Pi \ra$, it follows that there is some $\alpha \in \Pi$ such that
$\varphi(U_\alpha)$ is nontrivial. Let $u \in U_\alpha$ be such that $\varphi(u) \neq 1$. Since $W$ is infinite
and $(W,S)$ is irreducible, there exists $\beta \in \Phi$ such that $\alpha \cap \beta = \varnothing$
\cite[Proposition 8.1 p.309]{HDRHee}. Let $t = r_\beta r_\alpha \in W$ and $\alpha_n = t^n(\alpha)$ for all $n
\geq 0$. By definition, we have $\displaystyle \lim_{n \to +\infty} d(c_+, -\alpha_n) = +\infty$. Let $\tau \in
N$ be such that $\nu(\tau) = t \in W$, where $\nu : N \to N/T = W$ is the canonical projection. For each $n \geq
0$, let $u_n = \nu^n.u.\nu^{-n}$. Since $G$ is abelian, we have $\varphi(u_n) = \varphi(u) \neq 1$ for all $n$.
On the other hand, by definition $u_n \in U_{\alpha_n}$ and, hence $\displaystyle \lim_{n \to + \infty} u_n = 1$
by (FPRS). This contradicts the continuity of $\varphi$.
\end{proof}

The following lemma will be used again below, in order to establish restrictions on finite quotients of a group
endowed with a twin root datum.

\begin{lem}\label{lem:Rank1}
Let $(X, \{U_\alpha, U_{-\alpha}\})$ be a twin root datum of rank one. We have the following:
\begin{itemize}
\item[(i)] The group $X$ is not nilpotent.

\item[(ii)] Given a homomorphism $\varphi : X \to G$ whose kernel does not centralize $X_\alpha = \la U_\alpha
\cup U_{-\alpha} \ra$, we have $\varphi(U_\alpha) = \varphi(U_{-\alpha})$.
\end{itemize}
\end{lem}
\begin{proof}
(i). The  building associated with a twin root datum of rank one is merely an abstract set, which may be
identified with the conjugacy class of $U_\alpha$ in $X$. Clearly the kernel of the $X$-action is the
centralizer $Z_X(X_\alpha)$. This implies that the group $X_\alpha$ is not nilpotent: up to dividing $X_\alpha$
by its center, we obtain a group endowed with a twin root datum of rank one which acts faithfully on the
associated building, and is therefore center-free.

\vspace{2mm}

(ii). Since $(X, \{U_\alpha, U_{-\alpha}\})$ is a twin root datum of rank one, it follows that $X$ acts
$2$--transitively, hence primitively, on the conjugacy class of $U_\alpha$ in $X$. Therefore, the normal
subgroup $\Ker(\varphi)$ is transitive on the conjugacy class of $U_\alpha$ in $X_\alpha$ because it does not
centralize $U_\alpha$. In particular, this proves that $\varphi(U_\alpha) = \varphi(U_{-\alpha})$.
\end{proof}

\subsection{Topological simplicity} The following proposition is an improvement of the topological
simplicity theorem of \cite{Rem04} (see also \cite[Theorem~3.2]{CER}). We also note that, under some additional
assumptions, topological completions of Kac-Moody lattices have recently been shown to be \emph{abstractly}
simple by L.~Carbone, M.~Ershov and G.~Ritter \cite{CER}.

\begin{prop}
\label{prop:TopoSimplicity} Let $(W,S)$ be an irreducible Coxeter system of non-spherical type with associated
root system $\Phi$. Let $(\Lambda, \{U_\alpha\}_{\alpha \in \Phi})$ be a twin root datum of type $(W,S)$ with
finite root groups and let $\overline\Lambda_+$ be its positive topological completion. We assume that the root
groups are all solvable and that $[\overline\Lambda_+,\overline\Lambda_+]$ is dense in $\overline{\Lambda}_+$.
Then:
\begin{enumerate}
\item[(i)] Every closed subgroup of $\overline \Lambda_+$ normalized by $\Lambda^\dagger$ either contains
$\Lambda^\dagger$ or centralizes $\Lambda^\dagger$. In particular, the group
$\overline\Lambda^\dagger_+/Z(\Lambda^\dagger)$ is topologically simple.

\item[(ii)] Let $J$ be an irreducible non-spherical type in $S$ and let  $\overline G_J$ be the closure in
$\overline\Lambda_+$ of the group  generated by the root groups indexed by the simple roots in $J$ and their
opposites. Assume that $[\overline G_J,\overline G_J]$ is dense in $\overline G_J$. Then any proper closed
normal subgroup of $\overline G_J$ is contained in the center $Z(\overline G_J)$.
\end{enumerate}
\end{prop}

\begin{remark}
This is the opportunity to correct a mistake in \cite[Proposition 2.B.1  (iv)]{Rem04}. The factor groups there
are not topologically simple but simply have  property (ii) above: their proper closed normal subgroups fix
inessential buildings, but this does not seem to imply easily that the  whole ambient building is fixed. This
does not affect the rest of the paper. The second author thanks M.~Ershov for pointing out this mistake to him.
\end{remark}

\begin{proof}
For both (i) and (ii), the proof is an easy ``topological'' adaptation of the ``abstract'' arguments of
\cite[IV.2.7]{BbkLie4-5-6}. The essential point is that the group $\overline U_+$, and hence also $\overline G_J
\cap \overline U_+$, is pro-solvable by Proposition~\ref{prop:pronilpotency}(ii).
\end{proof}

\vspace{1cm}

\section{Non-affine Coxeter groups}
\label{s - non-affine Coxeter} This section is mainly Coxeter theoretic. We prove that in any non-affine
infinite Coxeter complex, given any  root there exist two other roots such that any two roots in the
so-obtained triple have empty intersection. Such a triple is called a  \emph{fundamental hyperbolic
configuration} and used in the next section to prove strong restrictions on finite index normal  subgroups for
twin root data.

\subsection{Parabolic closure}
\label{ss - parabolic closure} Let $(W,S)$ be a Coxeter system. Given a subset $R$ of $W$, we denote by $\Pc(R)$
the \textit{parabolic  closure} of $R$, namely the intersection of all parabolic subgroups of  $W$ containing
$R$. This notion is defined in D.~Krammer's PhD \cite{Kra94}. It is itself a parabolic subgroup which can be
characterized  geometrically as follows. Let $\mathscr{C}$ be the Coxeter complex associated with $(W,S)$. Given
$R \subset W$ and any simplex $\rho$ of maximal dimension stabilized by  $\langle R \rangle$, we have:
$\Pc(R)=\Stab_W(\rho)$.

\vspace{2mm}

By a \emph{Euclidean triangle group}, we mean a reflection subgroup of $\Isom(\E^2)$ which is the automorphism
group of a regular tessellation of the Euclidean plane $\E^2$ by triangles. Recall that there are three
isomorphism classes of such groups, corresponding respectively to tessellations by triangles with angles
$(\frac{\pi}{3}, \frac{\pi}{3}, \frac{\pi}{3})$, $(\frac{\pi}{2}, \frac{\pi}{4}, \frac{\pi}{4})$,
$(\frac{\pi}{2}, \frac{\pi}{3}, \frac{\pi}{6})$.

\begin{lem}\label{lem:Coxeter}
Let $(W,S)$ be a Coxeter system and let $r, s$ be reflections in $W$.
Assume that the product $\tau = rs$ is of infinite order.
Then  the following holds.
\begin{itemize}
\item[(i)] The Coxeter diagram of $\Pc(\tau)$ is irreducible.

\item[(ii)] The reflections $r$ and $s$ belong to $\Pc(\tau)$.

\item[(iii)] Let $t$ be a reflection which does not centralize $\tau$  and such that $\langle r,s,t \rangle$ is
isomorphic to a Euclidean  triangle group. Then $t$ belongs to $\Pc(\tau)$.
\end{itemize}
\end{lem}

\begin{proof}
We first prove that (ii) implies (i). Let us assume that (ii) holds. By a suitable conjugation in $W$, we may --
and shall -- assume that  $\Pc(\tau)=W_J$ for some subset $J \subseteq S$. Let $J_r$ be the connected component
of $J$ such that the irreducible  factor $W_{J_r}$ contains $r$. If $s$ did not belong to $W_{J_r}$ then $r$ and
$s$ would generate a  subgroup isomorphic to ${\bf Z}/2{\bf Z} \times {\bf Z}/2{\bf Z}$,  contradicting that
$\tau$ is of infinite order. Therefore $s \in  J_r$ and by definition of the parabolic closure we  have $J=J_r$.

\vspace{2mm}

Suppose (ii) fails. Without loss of generality, this means that $r \not \in \Pc(\tau)$. Let $\mathscr{C}$ be the
Coxeter complex associated with $(W,S)$ and $\rho$  be a simplex of maximal dimension which is stabilized by
$\langle\tau\rangle$. Since $r \not \in \Pc(\tau)$, it follows that $\rho$ is contained in the interior of one
of the two half-spaces determined by $r$. Let $\alpha$ be this half-space. We have $\rho \subset \alpha$, and
hence $ \rho \subset \bigcap_{n \in  {\bf Z}} \tau^n.\alpha$ because $\rho$ is $\tau$-invariant. This is absurd
since $\bigcap_{n \in {\bf Z}} \tau^n.\alpha$ is empty.

\vspace{2mm}

The proof of (iii) is similar.
Let $t$ be a reflection which does not centralize $\tau$ and such that  $\langle r, s, t \rangle$ is isomorphic to a Euclidean triangle group.
Let $\beta$ be any of the two half-spaces associated with $t$.
Using the fact that $\tau$ does not centralize $t$, it is immediate to  check in the Euclidean plane that the intersection $\bigcap_{n \in {\bf  Z}} \tau^n.\beta$ is empty.
Hence the same argument as in the proof of (ii) can be applied and  yields $t \in \Pc (\tau)$.
\end{proof}

We also need the following result due to D.~Krammer.
It is a first evidence that non-affine Coxeter groups have some weak  hyperbolic properties.

\begin{prop}
\label{prop:Krammer} Let $(W,S)$ be an irreducible, non-affine Coxeter system. Let $w  \in  W$ be such that
$\Pc(w)=W$. Then the cyclic group generated by $w$ is of finite index in its  centralizer.
\end{prop}
\begin{proof}[Reference] This is \cite[Corollary~6.3.10]{Kra94}. \end{proof}

\begin{remark}
This result is of course false for affine Coxeter groups of rank $\ge  3$, since the centralizer of any
translation in such a group contains  the translation subgroup.
\end{remark}

\subsection{Fundamental hyperbolic configuration}
\label{ss - fundmental hyperbolic} The non-linearity proof in \cite[\S 4]{Rem04} makes crucial use (for a  very
specific case of Weyl groups) of the fundamental hyperbolic  configuration defined in the introduction of this
section. We prove here that the Coxeter complex of any infinite non-affine  irreducible Coxeter group contains
many such configurations. Note that an affine Coxeter complex does not contain any fundamental hyperbolic
configuration. We do not assume the generating set $S$ to be finite.

\begin{theorem}
\label{thm:Coxeter} Let $(W,S)$ be an irreducible non-affine Coxeter system and let $\mathscr{C}$ be the
associated Coxeter complex. Let $\alpha, \beta$ be two disjoint non-opposite root half-spaces of $\mathscr{C}$.
Then there exists a root half-space $\gamma$ such that $\gamma\cap\alpha = \gamma\cap\beta = \varnothing$.
\end{theorem}

\begin{proof}
Let us first deal with the case when $S$ is infinite. The pair $\{r_\alpha ; r_\beta \}$ is contained in  a
finitely generated standard parabolic subgroup of $W$: take explicit (minimal) writings of $r_\alpha$ and
$r_\beta$ in the generating system $S$; the union of all elements used in these writings defines a finite
non-spherical subdiagram. Up to adding a finite number of vertices to this subdiagram, we may assume that it is
irreducible and non-affine. The corresponding standard parabolic subgroup of $W$ is finitely generated and
contains $r_\alpha$ and $r_\beta$.

\vspace{2mm}

We henceforth assume that the generating system $S$ is finite and  denote by $\mid\! S \!\mid$ its cardinality.
We prove the assertion by induction on $\mid\! S \!\mid$. The roots $\alpha$ and $\beta$ being non-opposite, the
corresponding  reflections $r_\alpha$ and $r_\beta$ generate an infinite dihedral  subgroup in $W$. This
excludes $\mid\! S \!\mid \, =1$ and $\mid\! S \!\mid=2$, except  possibly when the two vertices are connected
by an edge labelled by  $\infty$. But since the latter diagram is affine, the induction starts at $\mid\! S
\!\mid \, =3$.

\vspace{2mm}

Assume first that $\mid\! S \!\mid \, =3$, i.e. that the Coxeter  diagram of $(W,S)$ is a triangle.
Denoting by $a$, $b$ and $c$ the labels of its edges, we have $a,b,c  \ge 3$, and also
${1\over a}+ {1\over b}+ {1\over c}<1$ because $(W,S)$ is non-affine.
Let $\mathbb{H}^2$ denote the hyperbolic plane and let $\mathscr{T}$ be  a geodesic triangle in $\mathbb{H}^2$ of angles ${\pi\over a}$,  ${\pi\over b}$ and ${\pi\over c}$ (an angle equal to 0 correspond to a  vertex in the boundary of $\mathbb{H}^2$).
It follows from Poincar\'e's polyhedron theorem that the reflection  group generated by $\mathscr{T}$ is isomorphic to $W$ and that the  so-obtained hyperbolic tiling is a geometric realization of the Coxeter  complex of $(W,S)$
\cite[IV.H]{Maskit}.
Thanks to this geometric realization, the result is then clear when  $\mid\! S \!\mid \, =3$.

\vspace{2mm}

Assume now that $\mid\! S \!\mid \, >3$ and that the result is proved  for any Coxeter system as in the theorem
and whose canonical set of  generators has less than $\mid\! S \!\mid$ elements. Denote by $\tau$ the infinite
order element $r_\alpha r_\beta$. Using a suitable conjugation, we may -- and shall -- assume that  $\Pc(\tau)$
is standard parabolic, i.e. $\Pc(\tau)= W_J$ for some $J  \subseteq S$. According to Lemma~\ref{lem:Coxeter},
the Coxeter system $(W_J, J)$ is  irreducible by (i) and we have $r_\alpha, r_\beta \in  W_J$ by  (ii). Then two
cases occur.

\vspace{2mm}

The first case is when $(W_J, J)$ is non-affine. By the induction hypothesis, we only have to deal with the case
$J= S$ and $W_J =W$. If any canonical generator in $S$ centralized $\tau$, then we would have  $W=Z_W(\tau)$;
but $\tau$ cannot be central in $W$ since $\tau = r_\alpha r_\beta$ does not centralize $r_\alpha$ and
$r_\beta$. Therefore there exists a reflection $t  \in  S$ such that $t$ does  not centralize $\tau$. Let $T$ be
the subgroup generated by $t$, $r_\alpha$ and $r_\beta$. If $T$ were isomorphic to a Euclidean triangle group,
then $Z_T(\tau)$  would contain a free abelian group of rank 2. This is impossible by
Proposition~\ref{prop:Krammer}. Therefore, $T$ is isomorphic to a hyperbolic triangle group and we can  conclude
as in the case $\mid\! S \!\mid \, =3$.

\vspace{2mm}

The remaining case is when $(W_J, J)$ is affine. Then $J$ is properly contained in $S$ because $W$ is non-affine
and there exists an element $s  \in S \setminus J$ which does not normalize $W_J$. In particular $s$ does not
centralize $\tau$ because $\Pc(\tau)=W_J$. Let $T'$ be the subgroup generated by $s$, $r_\alpha$ and $r_\beta$.
If $T'$ is isomorphic to a Euclidean triangle group, then  Lemma~\ref{lem:Coxeter} (iii) implies that $s  \in
W_J$, which is  excluded. Thus $T'$ is isomorphic to a hyperbolic triangle group and we are again  reduced to
the case $\mid\! S \!\mid \, =3$.
\end{proof}

\section{Simplicity of twin building lattices}
\label{s - simple lattices}

As mentioned in the introduction, the proof of the main simplicity theorem applies to the general setting of
twin building lattices: the only required assumption is that root groups are nilpotent (see
Theorem~\ref{thm:SimpleDRJ}). The proof splits into two parts, each of which is presented in a separate
subsection below. These two parts have each their own specific hypotheses and are each of independent interest.

\subsection{Finite quotients of groups with a twin root datum}
\label{ss - no finite quotient} Here, we prove strong restrictions on finite index normal subgroups of a group
endowed with a twin root datum, under the assumption that root groups are nilpotent. These conditions are
fulfilled by Kac-Moody groups over arbitrary fields since their Levi factors are abstractly isomorphic to
reductive algebraic groups.

\begin{theorem}
\label{thm:NoFiniteQuotient} Let $(W,S)$ be a Coxeter system with associated root system $\Phi$ and let $\Pi$ be
the root basis associated to $S$. Let $G$ be a group endowed with a twin root datum $\{U_\alpha \}_{\alpha \in
\Phi}$ indexed by $\Phi$. Suppose that:
\begin{itemize}
\item[(1)] The Coxeter system $(W,S)$ is irreducible, non-spherical and non-affine;

\item[(2)] For any $\alpha \in \Pi$, the root group $U_\alpha$ is  nilpotent.
\end{itemize}
Let $H$ be a normal subgroup of $G$ such that $N/T.(N \cap H)$ is  finite. Let $G^\dagger$ be the subgroup of
$G$ generated by the root groups, let $H^\dagger = H \cap G^\dagger$, let $\pi : G^\dagger \to
G^\dagger/H^\dagger$ be the canonical projection and for each $\alpha \in \Pi$, denote by $f_\alpha$ the
inclusion $U_\alpha \to G^\dagger$. Then the composed map:

\vspace{2mm}

\centerline{$\begin{CD}
\prod_{\alpha \in \Pi} U_\alpha @>> \ \prod f_\alpha \ > G^\dagger @>> \pi > G^\dagger/H^\dagger
\end{CD}$}

\vspace{2mm}

is a surjective homomorphism. In particular, the group $G^\dagger/H^\dagger$ is nilpotent.
\end{theorem}

\begin{remark}
The finiteness of $N/T.(N \cap H)$ is automatically satisfied  when $H$ has finite index in $G$.
\end{remark}

\begin{figure}
\input{NoFiniteQuotient.pstex_t}
\caption{Proof of Theorem~\ref{thm:NoFiniteQuotient}}\label{Figure}
\end{figure}
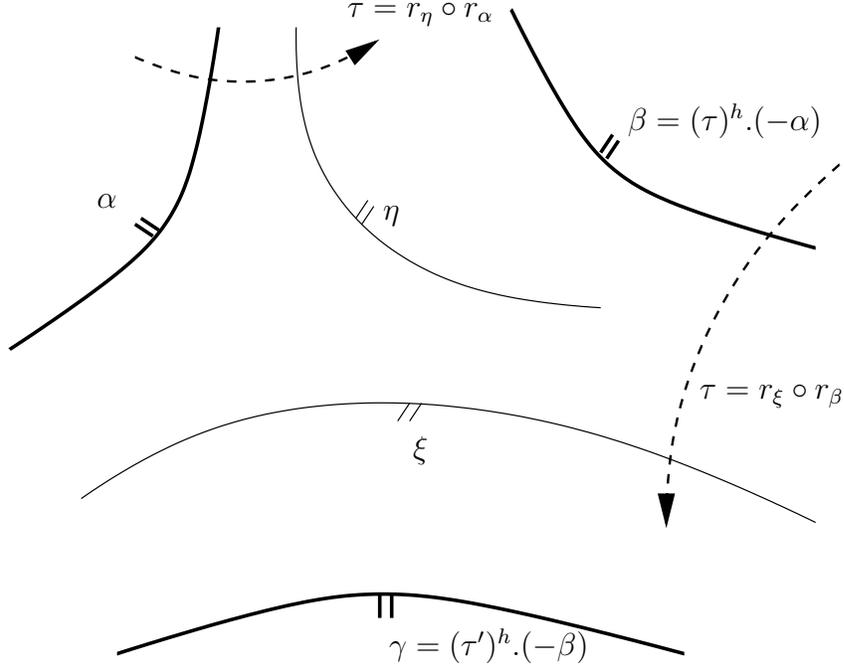

\begin{proof}
We identify the elements of $\Phi$ with the half-spaces of the  Davis complex $\mathscr{A}$ associated with
$(W,S)$. We set $h=[N:T.(N \cap H)]$.

\vspace{2mm}

Let $\alpha$ be an arbitrary root. By \cite[Proposition 8.1  p.309]{HDRHee} there is a root $\eta$ such that
$\alpha\cap\eta = \varnothing$. The product $\tau = r_\eta r_\alpha$ has infinite order. We set $\beta =
\tau^h.(-\alpha)  \in  \Phi$. We have $\beta \subset \eta$ and, hence, the roots $\alpha$ and  $\beta$ are
disjoint (see Figure on p.~\pageref{Figure}). By Theorem~\ref{thm:Coxeter}, there exists a root  $\xi  \in \Phi$
such that $\alpha\cap\xi = \eta \cap\xi = \varnothing$. In particular $\beta \cap \xi = \varnothing$. Again the
product $\tau'=r_\xi r_\beta$ has infinite order. We set $\gamma = (\tau')^h.(-\beta)$.

\vspace{2mm}

By construction, we have $\gamma \subset \xi$ (see Figure~\ref{Figure}). Hence the roots $\alpha$, $\beta$ and
$\gamma$ are pairwise disjoint. Therefore it follows from Assumption~(2) and
Proposition~\ref{prop:pronilpotency}(iii) that the group $U'=\langle U_\alpha \cup U_{-\gamma}\rangle$ is
nilpotent, and so is its image $\pi(U')$. But by (TRD2) we have

\vspace{2mm}

\centerline{$U_\beta = \tau^h U_{-\alpha} \tau^{-h}
\hspace{.5cm}  \text{ and } \hspace{.5cm} U_{-\gamma}  = (\tau')^h  U_\beta (\tau')^{-h}$.}

\vspace{2mm}

Note that $N/T.(N \cap H)$ is the quotient of the Weyl group $W= N/T$  induced by $H \triangleleft G$. Since $h$
is the order of the quotient $N/T.(N \cap H)$, applying $\pi$  provides $\pi(U_{-\alpha})= \pi(U_\beta)=
\pi(U_{-\gamma})$, which  implies

\vspace{2mm}

\centerline{$\pi(U') = \langle\pi(U_\alpha) \cup \pi(U_{-\alpha})  \rangle= \pi(X_\alpha)$.}

\vspace{2mm}

This shows that $\pi(X_\alpha)$  is a nilpotent group.

\vspace{2mm}

Note that $(X_\alpha, \{U_\alpha, U_{-\alpha}\})$ is a twin root datum of rank one. By Lemma~\ref{lem:Rank1}(i),
the group $X_\alpha$ is not nilpotent and, hence, $X_\alpha\cap H={\rm Ker}(\pi \!\mid_{X_\alpha})$ is not
central in $X_\alpha$. Therefore, we have $\pi(U_\alpha) = \pi(U_{-\alpha})$ by Lemma~\ref{lem:Rank1}(ii).

\vspace{2mm}

Finally, for any two distinct roots $\alpha, \beta \in \Pi$, we have  $[U_\alpha, U_{-\beta}]=1$ by axiom
(TRD1). In view of the preceding paragraph, this implies that $[\pi(U_\alpha), \pi(U_{\beta})]=1$ for all
distinct $\alpha, \beta \in \Pi$. The desired result follows by noticing that $G^\dagger$ is generated by
$\bigcup_{\pm \alpha \in \Pi} U_\alpha$. This is easily seen using axiom (TRD2) of twin root data to produce
elements in $N$ and then to  conjugate the simple root groups by these elements to produce any desired root
group.
\end{proof}

\vspace{2mm}

\vspace{2mm}

The following corollary applies to all split and almost split Kac-Moody groups over finite fields.

\begin{cor}
\label{cor:NoFiniteQuotient} Let $G$ be a group as in Theorem \ref{thm:NoFiniteQuotient}, maintain the
assumptions (1) and (2) and assume moreover that root groups are finite. Here we let $H^{\dagger}$ denote the
intersection of all finite index normal subgroups of $G^\dagger$. Then

\vspace{2mm}

\centerline{$[G^\dagger : H^{\dagger}] \, \leq \, \prod_{\alpha\in \Pi} |U_\alpha|$.}

\vspace{2mm}

Furthermore, we have $H^{\dagger} = G^\dagger$ whenever one of the following holds:
\begin{itemize}
\item[(3)] Each group $X_{\alpha}$, $\alpha  \in  \Pi$, is a finite group of Lie type and the minimal order
$q_{\rm min} = \min \{ |U_\alpha| : \alpha   \in  \Pi \} > 3$;

\item[(4)] The Coxeter system $(W,S)$ is $2$-spherical, i.e. every $2$--subset of $S$ generates a finite group,
and $q_{\rm min} > 2$;

\item[(5)] The Coxeter system $(W,S)$ is simply laced, i.e. every $2$--subset of $S$ generates a group of
order~$4$ or~$6$.
\end{itemize}
\end{cor}

\begin{remark}
The above group $H^{\dagger}$ is contained in any finite index  normal subgroup of $G$.
\end{remark}

\begin{proof}
Let $H$ be a finite index normal subgroup in $G^\dagger$. Applying Theorem \ref{thm:NoFiniteQuotient} to
$G^\dagger$ we see that the index $[G^\dagger:H]$ is uniformly bounded, so that finite index subgroups of
$G^\dagger$ are finite in number. This implies that the intersection defining $H^\dagger$ is finite, so that
$H^\dagger$ is itself a finite index subgroup. It remains to apply again Theorem~\ref{thm:NoFiniteQuotient} to
obtain the desired bound on $[G^\dagger : H^{\dagger}]$.

\vspace{2mm}

If condition (3) holds, then each $X_\alpha$ is perfect (in fact: simple modulo center), so admits no
non-trivial nilpotent quotient. Equality $H^{\dagger}= G^\dagger$ follows from
Theorem~\ref{thm:NoFiniteQuotient} applied to $G^\dagger$. Similarly, if (4) or (5) holds then for each $\alpha
\in \Pi$ there  exists $\beta \in \Pi- \{\alpha\}$ such that $X_{\alpha, \beta}=\la X_\alpha, X_\beta \ra$ is a
rank~$2$ finite  group of Lie type. All such groups are perfect except $B_2(2)$ and $G_2(2)$ (which contain both
a simple subgroup of index~$2$). Since (4) implies $q_{\rm min}
> 2$ and (5) implies that $X_{\alpha, \beta}$ is of type $A_2$, the group $X_{\alpha, \beta}$ is isomorphic  to
neither of the latter groups and we have again $H^{\dagger} = G^\dagger$.
\end{proof}

Theorem~\ref{thm:NoFiniteQuotient} and its corollary imply that any Kac-Moody group over a finite field of
irreducible non-spherical and non-affine type, admits at most a finite  number of finite quotients, which are
necessarily abelian. Furthermore, if the ground field is of cardinality at least~$4$ and if the group is
generated by its root groups, e.g. because it is simply  connected, then all finite quotients are trivial.

\vspace{2mm}

We close this subsection with an example of a Kac-Moody group which  admits nontrivial finite quotients when the
ground field is ${\bf  F}_2$ or ${\bf F}_3$. We set $I=\{1, 2, 3\}$ and consider the generalized Cartan matrix

\vspace{2mm}

\centerline{$A = (A_{ij})_{i,j \in I} = \left(\begin{array}{ccc}
\hfill 2 & -2 & -2\\
-2 & \hfill 2 & -2\\
-2 & -2 & \hfill 2
\end{array}\right)$.}

\vspace{2mm}

Let $\mathcal{G}_A$ be the simply connected Tits functor of type $A$  \cite[3.7.c]{TitsJA}. We set $\Lambda =
\mathcal{G}_A({\bf F}_2)$. For each $i  \in  I$, we let $\varphi_i : {\rm SL}_2({\bf F}_2) \to  \Lambda$ be the
standard homomorphism \cite[\S 2 and 3.9]{TitsJA} and  let $f_i : {\rm SL}_2({\bf F}_2) \to {\bf F}_2$ be the
surjective  homomorphism defined by $f_i\left(\begin{array}{cc}
1 & 1 \\
0 & 1
\end{array}\right)
= f_i\left(\begin{array}{cc}
1 & 0 \\
1 & 1
\end{array}\right) = 1$.
Using the defining relations of $\Lambda$ \cite[\S 8.3]{RemAst}, we see  that there is a unique homomorphism $f
: \Lambda \to \prod_{i \in I} {\bf F}_2$ such that $f \circ  (\prod_{i \in I} \varphi_i)=\prod_{i \in I} f_i$.
By the definition of  $f_i$, the homomorphism $f$ is surjective.

\subsection{Non-arithmeticity}

For the next statement we recall that for a group inclusion $A<B$, the  \emph{commensurator} of $A$ in $B$,
denoted ${\rm Comm}_B(A)$, consists  of the elements $b  \in  B$ such that $A$ and $bAb\inv$ share a  finite
index subgroup. According to a well-known theorem of G.~Margulis, a lattice in a semisimple Lie group is
arithmetic if and only if its commensurator is dense in the ambient Lie group \cite[6.2.5]{Zimmer}.

\begin{cor}
\label{cor - non arithmetic} Let $\Lambda$ be a group as in Theorem \ref{thm:NoFiniteQuotient}. Suppose moreover
that assumption (3) of Corollary~\ref{cor:NoFiniteQuotient} holds and let $\overline \Lambda_+$ (resp.
$\overline \Lambda_-$) be the positive (resp. negative) topological completion of $\Lambda$. Then the
commensurator ${\rm Comm}_{\overline\Lambda_+ \times \overline\Lambda_-}\bigl(\Lambda^\dagger \bigr)$ is a
discrete subgroup of $\overline\Lambda_+ \times \overline\Lambda_-$.
\end{cor}

\begin{proof}
Recall that $\Lambda^\dagger$ is viewed here as a diagonal subgroup of $\overline\Lambda_+ \times
\overline\Lambda_-$. By Corollary~\ref{cor:NoFiniteQuotient}, the commensurator ${\rm Comm}_{\overline\Lambda_+
\times \overline\Lambda_-}\bigl(\Lambda^\dagger \bigr)$ is equal to the normalizer $N_{\overline\Lambda_+ \times
\overline\Lambda_-}\bigl(\Lambda^\dagger \bigr)$ because any finite index subgroup of a given group contains a
finite index normal subgroup. Furthermore, the centralizer $Z_{\overline \Lambda_+}(\Lambda^\dagger)$ (resp.
$Z_{\overline \Lambda_-}(\Lambda^\dagger)$) is nothing but the kernel of the $\overline \Lambda_+$-action (resp.
$\overline \Lambda_-$-action) on the positive (resp. negative) building associated with $\Lambda$. By
Proposition~\ref{prop:TopoCompletions}(iii), we have $Z_{\overline \Lambda_+}(\Lambda^\dagger)=
Z_\Lambda(\Lambda^\dagger)$ (resp. $Z_{\overline \Lambda_-}(\Lambda^\dagger)= Z_\Lambda(\Lambda^\dagger)$).
Therefore, we have an exact sequence:

\vspace{2mm}

\centerline{$1 \longrightarrow Z_\Lambda(\Lambda^\dagger) \times Z_\Lambda(\Lambda^\dagger) \longrightarrow
N_{\overline\Lambda_+ \times \overline\Lambda_-}\bigl(\Lambda^\dagger \bigr) \longrightarrow
\Aut(\Lambda^\dagger)$.}

\vspace{2mm}

This yields an exact sequence

\vspace{2mm}

\centerline{$1 \longrightarrow \bigl(Z_\Lambda(\Lambda^\dagger) \times Z_\Lambda(\Lambda^\dagger)\bigr).\Lambda^\dagger
\longrightarrow
N_{\overline\Lambda_+ \times \overline\Lambda_-}\bigl(\Lambda^\dagger \bigr) \longrightarrow
\Out(\Lambda^\dagger)$,}

\vspace{2mm}

where $\Out(\Lambda^\dagger) = \Aut(\Lambda^\dagger)/\Inn(\Lambda^\dagger)$ is the outer automorphism group. By
\cite[Corollary~B]{CM04}, the group $\Out(\Lambda^\dagger)$ is finite. Since $(Z_\Lambda(\Lambda^\dagger) \times
Z_\Lambda(\Lambda^\dagger)).\Lambda^\dagger$ is a discrete subgroup of ${\overline\Lambda_+ \times
\overline\Lambda_-}$ by Propositions~\ref{prop:TopoCompletions}(iii) and~\ref{prop:lattices}, it finally follows
that ${\rm Comm}_{\overline\Lambda_+ \times \overline\Lambda_-}\bigl(\Lambda^\dagger \bigr) =
N_{\overline\Lambda_+ \times \overline\Lambda_-}\bigl(\Lambda^\dagger \bigr)$ is discrete as well.
\end{proof}

\subsection{Normal subgroup property}
\label{ss - NSP} In view of Corollary~\ref{cor:NoFiniteQuotient}, the complementary property necessary in order
to obtain simplicity of the group $H^{\dagger}$ (modulo center) is that any non-central normal subgroup has
finite index. This is called the \emph{normal subgroup property} and is well-known for irreducible higher rank
lattices in Lie groups \cite[IV.4.9]{Margulis}. The generalization to irreducible cocompact lattices in products
of topological groups follows from work by U. Bader and Y. Shalom, following Margulis' general strategy (see
\cite{Shalom} and mostly \cite[Introduction]{BSh06} for an explanation of the substantial differences with the
classical case). In the attempt of adapting these results to Kac-Moody groups over finite fields, one has to
overcome the fact that Kac-Moody lattices are never cocompact. This was done in \cite{Rem05} by proving that one
can find a fundamental domain $D$ for $\Lambda$ in $\overline\Lambda_+ \times \overline\Lambda_-$ with respect
to which the associated cocycle is square integrable.

\vspace{2mm}

The following result is a restatement of the normal subgroup theorem proved in  \cite{BSh06} and \cite{Rem05},
in the general framework of twin building lattices.

\begin{theorem}
\label{th - NSP} Let $(\Lambda, \{U_\alpha\}_{\alpha \in \Phi})$ be a twin root datum of type $(W,S)$. We assume
that:
\begin{description}
\item[(NSP1)] Each root group $U_\alpha$ is finite and solvable.

\item[(NSP2)] The group $\Lambda$ is a twin building lattice in the sense of \ref{sect:TBLattice}.

\item[(NSP3)] The derived group $[\overline\Lambda_+, \overline\Lambda_+]$ is dense in $\overline \Lambda_+$.
\end{description}
Then any subgroup of $\Lambda$ which is normalized by $\Lambda^\dagger$, either centralizes $\Lambda^\dagger$ or
contains a finite index subgroup of $\Lambda^\dagger$.
\end{theorem}

\begin{remarks}
1. There is no condition excluding affine diagrams. Indeed, Kac-Moody groups of affine type are
$\{0,\infty\}$-arithmetic groups and as such are irreducible lattices in higher-rank algebraic groups:  this
case was already covered by Margulis' theorem.

2. As pointed out to us by M.~Burger, an infinite group with the normal subgroup property cannot be hyperbolic since it is incompatible with  SQ-universality, the property that any countable group embeds in a suitable quotient of
the group under consideration. Any non-elementary hyperbolic group is SQ-universal
\cite{Delzant}, \cite{Olshanskii}. The fact that no Kac-Moody group can be hyperbolic can also be derived from the specific
property that Kac-Moody groups over finite fields contain infinitely many conjugacy classes of torsion elements.
\end{remarks}

\begin{proof}
Let $H$ be a normal subgroup of $\Lambda$ which does not centralize $\Lambda^\dagger$ and set $H^\dagger = H
\cap \Lambda^\dagger$. Note that a subgroup of $\Lambda$ (resp. $\overline \Lambda_+$) centralizes
$\Lambda^\dagger$ if and only if its acts trivially on the building $\mathscr{B}_+$. Hence, it suffices to show
that the index of $H^\dagger $ in $ \Lambda^\dagger$ is finite. To this end, we apply the main results of
\cite{BSh06}. This requires to ensure that two conditions are fulfilled. The first condition is that the closure
of $H^\dagger$ in $\overline \Lambda_\pm^\dagger$ is cocompact; this is an immediate consequence of
Proposition~\ref{prop:TopoSimplicity}(i). The second condition is the existence of a fundamental domain $D$ for
$\Lambda$ with respect to which the associated cocycle is square integrable; this is provided by the same
arguments as in \cite{Rem05}. We do not go into details here because this question is more carefully examined
Subsect.~\ref{ss - uniform integrability}, where we prove a refinement of the square integrability. We merely
remark that the group combinatorics needed to prove the existence of $D$, namely the structure of refined Tits
system defined in \cite{KacPet}, is available for arbitrary twin root data, and not only for those arising from
Kac-Moody groups, see Proposition~\ref{prop:TopoCompletions}(vi).
\end{proof}

\subsection{Simplicity of lattices}
\label{ss - simple lattices} We can now put together the two ingredients needed to prove the  simplicity theorem
for twin building lattices.

\begin{theorem}
\label{thm:SimpleDRJ} Let $(\Lambda, \{U_\alpha\}_{\alpha \in \Phi})$ be a twin root datum of type $(W,S)$. Let
$\Lambda^\dagger$ be the subgroup generated by all root groups and assume that:
\begin{description}
\item[(S0)] The Coxeter system $(W,S)$ is irreducible, non-spherical and non-affine.

\item[(S1)] Each root group $U_\alpha$ is finite and nilpotent.

\item[(S2)] The group $\Lambda$ is a twin building lattice in the sense of \ref{sect:TBLattice}.

\item[(S3)] The derived group $[\overline\Lambda_+, \overline\Lambda_+]$ is dense in $\overline \Lambda_+$.
\end{description}
Then the quotient $\Lambda^\dagger/Z(\Lambda^\dagger)$ is virtually simple.

\vspace{2mm}

Assume instead of (S3) that:
\begin{description}
\item[(S3+)] Each rank one group $X_\alpha= \la U_\alpha \cup U_{-\alpha} \ra$  is perfect.
\end{description}
Then any subgroup of $\Lambda$, normalized by $\Lambda^\dagger$, either centralizes $\Lambda^\dagger$ or
contains $\Lambda^\dagger$.
\end{theorem}

\begin{remark}
This theorem applies to the groups of mixed characteristics defined in \cite{RR06} provided the minimal size of
the finite ground fields is large enough with respect to the growth of the (right-angled) Weyl group. In this
case the lattices are by definition generated by their root groups and condition (S3+) is fulfilled. One can
push a little further this construction by replacing the rank~$1$ Levi factors, isomorphic to some suitable
${\rm SL}_2(q)$'s, by affine groups. In this case, the root groups are isomorphic to multiplicative groups of
finite fields, so the thicknesses are prime powers, and rank~$1$ subgroups are solvable. Condition (S3) is still
fulfilled in view of Proposition~\ref{prop:FPRS} and Lemma~\ref{lem:density}.
\end{remark}

\begin{proof}
Let $H^{\dagger}$ be the intersection of all finite index normal subgroups of $\Lambda^\dagger$. The center of
$H^{\dagger}$ is a normal subgroup of $\Lambda^\dagger$, which must be central in $\Lambda^\dagger$ in view of
Theorem~\ref{th - NSP}. In particular $Z(H^{\dagger})$ is finite and, moreover, the canonical projection of
$H^\dagger$ in $\Lambda^\dagger/Z(\Lambda^\dagger)$ is isomorphic to $H^\dagger/Z(H^\dagger)$ and coincides with
the intersection of all finite index subgroups of $\Lambda^\dagger/Z(\Lambda^\dagger)$.

\vspace{2mm}

By Corollary~\ref{cor:NoFiniteQuotient}, the index of $H^{\dagger}$ in $\Lambda^\dagger$ is finite. On the other
hand, it follows from Theorem~\ref{th - NSP} that $\Lambda^\dagger/Z(\Lambda^\dagger)$ is just infinite (i.e.
every nontrivial quotient is finite). Therefore, it follows from \cite[Proposition 1]{Wilson} that
$H^\dagger/Z(H^\dagger)$ is a direct product of finitely many isomorphic simple groups. Write
$H^\dagger/Z(H^\dagger) = H_1 \times \dots \times H_k$. We must prove that $k=1$.

\vspace{2mm}

Notice that $H^\dagger$, viewed as a diagonally embedded subgroup, is a lattice in $\overline \Lambda^\dagger_+
\times \overline \Lambda^\dagger_-$, because it is a finite index subgroup of the lattice $\Lambda^\dagger$.
Furthermore, $H^\dagger$ is irreducible. Indeed, since $H^\dagger$ is a finite index normal subgroup of
$\Lambda^\dagger$, its closure $\overline H^\dagger$ in $\overline \Lambda_+^\dagger$ is a non-central closed
normal subgroup, which must coincide with $\overline \Lambda_+^\dagger$ by
Proposition~\ref{prop:TopoSimplicity}(i).

\vspace{2mm}

Assume now that $k>1$. It follows that the simple group $H_1$ is a quotient of $H^\dagger$ which is not
co-central,
since we have a composed map
$$H^\dagger \to H^\dagger/Z(H^\dagger)  = H_1 \times \dots \times H_k \to H_1.$$
The closure of the projection of the corresponding normal subgroup of $H^\dagger$ in $\overline
\Lambda_+^\dagger$ is thus a non-central closed normal subgroup of $\overline \Lambda_+^\dagger$. Hence it
coincides with $\overline \Lambda_+^\dagger$ by Proposition~\ref{prop:TopoSimplicity}(i). By
\cite[Theorem~1.3]{BSh06}, this implies that $H_1$ is amenable. Since the $H_i$'s are all isomorphic, it follows
that $H^\dagger/Z(H^\dagger)$ is amenable, and so is $\Lambda^\dagger$ since $Z(H^\dagger)$ and
$[\Lambda^\dagger: H^\dagger]$ are finite. Recall that $\Lambda^\dagger$ acts on the associated positive
building $\mathscr{B}_+$, which may be viewed as a proper ${\rm CAT}(0)$-space. Amenability of $\Lambda^\dagger$
implies that its action on $\mathscr{B}_+$ stabilizes a Euclidean flat or fixes a point in the visual boundary
at infinity \cite{AB98}. Both eventualities are absurd. This shows that $k=1$ as desired.

\vspace{2mm}

Assume now that (S3+) also holds. Note that a subgroup of $\Lambda$ (resp. $\overline \Lambda_+$) centralizes
$\Lambda^\dagger$ if and only if its acts trivially on the building $\mathscr{B}_+$. Hence, in view of what has
already been proven, it suffices to show that $H^\dagger = \Lambda^\dagger$. This follows from
Corollary~\ref{cor:NoFiniteQuotient}.
\end{proof}

Here is now the Kac-Moody specialization of this theorem:

\begin{theorem}
\label{th:Simple lattices} Let $\Lambda$ be a split or almost split Kac-Moody group over a finite  field ${\bf
F}_q$ of order $q$. Let us denote by $(W,S)$ the natural Coxeter system of the Weyl group  $W$ and by $W(t)$ the
growth series of $W$ with respect to $S$. Assume that $(W,S)$ is irreducible, neither of spherical nor of affine
type and that $W({1\over q})<+\infty$. Then the derived group of $\Lambda$, divided by its center, is simple.
\end{theorem}

\begin{proof}
All root groups of $\Lambda$ are nilpotent (of class at most~$2$). Thus conditions (S0), (S1) and (S2) are
clearly satisfied. In order to deduce the desired statement from Theorem~\ref{thm:SimpleDRJ} and its proof, it
remains to show that the derived group $[\Lambda, \Lambda]$ coincides with the intersection $H^\dagger$ of all
finite index subgroups of $\Lambda^\dagger$.

\vspace{2mm}

Each rank one subgroup $X_\alpha = \la U_\alpha \cup U_{-\alpha} \ra$ is isomorphic to the $\F_q$-points of a
simple algebraic group of relative rank one. Therefore, the group $X_\alpha$ is perfect except if $U_\alpha$ is
of order~$2$ or~$3$ in which case it is abelian. In view of Theorem~\ref{thm:NoFiniteQuotient}, this implies in
particular that the quotient $\Lambda / H^\dagger$ is abelian. Thus $[\Lambda, \Lambda] \subset H^\dagger$. It
follows from the proof of Theorem~\ref{thm:SimpleDRJ} that the latter inclusion cannot be proper, as desired.
\end{proof}

Note that if $q > 3$ then every  rank one subgroup of the Kac-Moody group $\Lambda$ is perfect and, hence,
condition (S3+) holds. In that case, we have $[\Lambda, \Lambda] = \Lambda^\dagger$.

\vspace{2mm}

Kac-Moody groups are the values over fields of group functors defined  by J.~Tits thanks to combinatorial data
called {\it Kac-Moody root  data} \cite{TitsJA}. The main information in a Kac-Moody root datum is given by a
generalized Cartan matrix, say $A$. Once $A$ is fixed we can still  make some choices in order for $\Lambda$  to
be generated by its root groups. In this case, e.g. when we choose the {\it simply connected Kac-Moody  root
datum}~\cite[3.7.c]{TitsJA},  we have  $\Lambda=\Lambda^\dagger=[\Lambda,\Lambda]$ (if $q>3$) and we recover the
situation described in the comment of the Simplicity theorem  (Introduction).

\vspace{2mm}

Let us now state a corollary on property~(T). Its proof is an opportunistic use of work by J.~Dymara and
T.~Januszkiewicz and by P.~Abramenko and B.~M\"uhlherr, but has nice  rigidity consequences (Theorem \ref{th -
action on CAT(-1)}).

\begin{cor}
\label{cor - T and simple} Let $\Lambda$ be a group endowed with a twin root datum satisfying (S0), (S1), (S2)
and (S3+) of Theorem \ref{thm:SimpleDRJ}. We assume furthermore that any two canonical reflections in $S$
generate a finite subgroup of $W$. If $q_{\rm min} > 1764^{|S|}$, then $\Lambda$ has Kazhdan's property (T). In
particular there exist infinitely many isomorphism classes of finitely presented infinite simple groups with
Kazhdan's property (T).
\end{cor}

\begin{proof}
This is a straightforward application of \cite[Theorem E]{DJ02}, which  provides the 1-cohomology vanishing
useful to a well-known criterion  for property (T) \cite[Chapitre 4]{dHV}. Finite presentation follows from
\cite{AbrMu} under the hypothesis that $q_{\min} > 2$. Finally, it follows from~\cite{CM04} that Kac-Moody
groups over non isomorphic finite fields (or of different types) are not isomorphic.
\end{proof}

Concretely, in order to produce infinite simple Kazhdan groups, it is  enough to pick a generalized Cartan
matrix $A=[A_{s,t}]_{s,t\in S}$ such that $A_{s,t}A_{t,s} \leq 3$ for each $s \neq t$ and a finite  ground
field, whose order is at least the size of $A$. The above simple groups seem to be the first examples of
infinite finitely generated simple groups enjoying property (T). The simple lattices in products of trees
constructed by M.~Burger and Sh.~Mozes \cite{BMProducts} are finitely presented but they cannot have  property
(T) since they act fixed point freely on trees.
However these lattices are torsion free, while a Kac-Moody group over a  finite field of characteristic $p$
contains infinite subgroups of exponent $p$ \cite[proof of Theorem 4.6]{RemNewton}.

\vspace{2mm}

\section{Non-linearity of Kac-Moody groups}
\label{s - non linear lattices}

A result of Mal'cev's asserts that any finitely generated linear group  is residually finite. In particular, the
groups covered by  Corollary~\ref{cor:NoFiniteQuotient} are not linear over any field. Note that with the
notation and assumptions of this corollary, the  group $G^\dagger$ is finitely generated. In this section, we
show that the latter corollary actually implies a  strong non-linearity statement for Kac-Moody groups over
arbitrary  fields of positive characteristic.

\subsection{Normal subgroups (arbitrary ground field)}
The proof of the non linearity theorem below  (Theorem~\ref{thm:NonLinearity}) requires the following statement,
which is a complement to Theorem~\ref{thm:NoFiniteQuotient}. The reader familiar with infinite-dimensional Lie
algebras will  recognize some similarity with \cite[Proposition 1.7]{Kac}.

\begin{prop}\label{prop:normal:InfiniteFields}
Let $A$ be a generalized Cartan matrix which is indecomposable and  non-affine, let $\mathcal{G}_A$ be a Tits
functor of type $A$ and let  $\K$ be an infinite field. We set $G = \mathcal{G}_A(\K)$ and $G^\dagger = \la
U_\alpha : \alpha   \in  \Phi \ra$, where $\{U_\alpha\}_{\alpha \in \Phi}$ is the twin  root datum given by the
root groups. Then given any normal subgroup $H$ of $G$, either $H$ contains  $G^\dagger$ or $H \cap U_\Psi =
\{1\}$ for each nilpotent set of roots  $\Psi$.
\end{prop}

\begin{proof}
Let $H \triangleleft G$ be such that $H \cap U_\Psi \neq \{1 \}$ for  some nilpotent set of roots $\Psi$. We
must prove that $H \supset G^\dagger$. The set $\Psi$ is finite \cite[2.2.6]{RemAst}; we assume that it is of
minimal cardinality with respect to the property that $H \cap U_\Psi\neq \{1 \}$ and we set $n=|\Psi|$.

\vspace{2mm}

Suppose that $n>1$. The elements of $\Psi$ can be ordered in a nibbling sequence $\alpha_1, \alpha_2, \dots,
\alpha_n$ [loc. cit., 1.4.1]. Now let $g  \in  H \cap U_\Psi - \{1 \}$. The group $U_\Psi$ decomposes as a
product $U_\Psi = U_{\alpha_1}  U_{\alpha_2} \dots U_{\alpha_n}$ [loc. cit., 1.5.2], so we have $g =  u_1 u_2
\dots u_n$ with $u_i \in U_{\alpha_i}$ for each $i=1$, ...,  $n$. By the minimality assumption on $\Psi$, the
elements $u_1$ and $u_n$  must be nontrivial. We set $j = \min \{i >1 :  u_i \neq 1\}$. By
Lemma~\ref{lem:nonaffinetorus} below, we can pick some $h  \in   Z_T(U_{\alpha_1})$ not centralizing
$U_{\alpha_j})$. By the defining relations of $\mathcal{G}_A$, we see in a suitable  parametrization of
$U_{\alpha_j}$ by the additive group $(\K,+)$ that  the action of $h$ on $U_{\alpha_j}$ by conjugation is merely
a  multiplication by an element of $\K^\times$. Therefore $h$ centralizes no nontrivial element of
$U_{\alpha_j}$ and  we obtain successively:

\vspace{2mm}

\centerline{$\begin{array}{rcl}
g\inv h\inv g h &=& u_n\inv \dots u_j \inv u_1\inv u_1^h u_j^h
\dots u_n^h\\
& = & u_n\inv \dots u_j\inv u_1\inv u_1 u_j^h \dots u_n^h\\
& = & u_n\inv \dots u_j\inv  u_j^h \dots u_n^h\\
& = & u_j\inv u_j^h u_{j+1}' \dots u_n'
\end{array}$}

\vspace{2mm}

for some $u_i' \in U_{\alpha_i}$ ($i = j+1, \dots, n$) and where the  last equality follows from the commutation
relations satisfied by the  $U_\alpha$'s in view of (TRD~1). Since $T$ normalizes $H$, we have $g\inv h\inv g h
\in  H$. Moreover the definition of $j$ and the choice of $h$ imply that $u_j  \inv u_j^h \neq 1$ so in
particular $g\inv h\inv g h \neq 1$. This shows that $H \cap U_{\Psi - \{\alpha_1\}} \neq \{1 \}$, which
contradicts the minimality of $|\Psi|$. Thus $n = 1$.

\vspace{2mm}

This shows that $H \cap U_\alpha$ is nontrivial for some $\alpha   \in  \Phi^+$. Since the group $\la U_\alpha
\cup U_{-\alpha}\ra$ is quasi-simple and  since its center intersects $U_\alpha$ trivially, we deduce that $H$
contains $U_\alpha \cup U_{-\alpha}$. Let $\Pi$ be a basis of $\Phi$ containing $\alpha$ and let $\beta \in  \Pi
-\{\alpha \}$ be such that the associated reflections $r_\alpha$  and $r_\beta$ do not commute. Since $\K$ is
infinite, it follows that $T \cap \la U_\alpha \cup  U_{-\alpha} \ra \not \subset Z_T(U_\beta)$. In particular,
there exists $h' \!\in T \cap H$ and $u  \in  U_\beta$  such that $ h'u (h')\inv u\inv \neq 1$. Thus $H \cap
U_\beta$ is nontrivial and, as above, this implies that  $H$ contains  $U_\beta \cup U_{-\beta}$. Finally, since
$A$ is indecomposable we obtain $U_\gamma < H$ for any  $\gamma  \in  \Phi$, that is to say $H \supset
G^\dagger$.
\end{proof}

Let us now immediately prove the lemma we used in the previous proof.

\begin{lem}
\label{lem:nonaffinetorus}
Maintain the notation and assumptions of  Proposition~\ref{prop:normal:InfiniteFields} and set
$T = \bigcap_{\alpha \in \Phi} N_G(U_\alpha)$.
Then, for any positive roots $\alpha, \beta \in \Phi^+$, the inclusion  $Z_T(U_\alpha) \subset Z_T(U_\beta)$ implies $\alpha = \beta$.
\end{lem}

\begin{proof}
Let $\Lambda$ (resp. $\Lambda^\vee$) be the lattice of algebraic  characters (resp. cocharacters) of the maximal
split torus $T$ \cite[\S 8.4.3]{RemAst}. We may -- and shall -- identify the abstract root system $\Phi$ to a
subset of $\Lambda$ and use the identification $T \simeq \Hom_{\rm groups}(\Lambda, \K^\times)$. For $\alpha \in
\Phi$, $x  \in  \K$ and $t  \in  T$, we have:
\renewcommand{\theequation}{$**$}
\begin{equation}\label{eq:*}
t.\mathfrak{u}_\alpha(x).t\inv = \mathfrak{u}_\alpha(t(\alpha).x),
\end{equation}

where $\mathfrak{u}_\alpha: (\K,+) \to U_\alpha$ is a standard isomorphism (see \cite[\S~3]{TitsJA}).

\vspace{2mm}

Assume now that $\alpha \neq \beta$. We claim that there exists $\gamma \in \Lambda^\vee$ such that $\la \alpha
\, |\, \gamma \ra = 0$ and $\la \beta \, |\, \gamma \ra \neq 0$. Let us set $k(\alpha, \beta)=\la \alpha \, |\,
\beta^\vee \ra \la \beta  \, |\, \alpha^\vee \ra$. If $k \neq 4$, then it is easy to see that  there exists such
a $\gamma$ in the group $\Z \alpha^\vee + \Z  \beta^\vee$. If $k=4$, then the order of $s_\alpha s_\beta$ is
infinite and we are  in position to apply Theorem~\ref{thm:Coxeter}. This yields a non-degenerate infinite
rank~$3$ root subsystem of $\Phi$  containing $\alpha$ and $\beta$. Then it is again easy to check the existence
of $\gamma$ inside this  subsystem. In both cases, the claim above holds. Given $t \in \K^\times$, let
$t^\gamma$ denote the element of $T$  defined by $t^\gamma : \lambda \mapsto t^{\la \lambda \, |\, \gamma \ra}$.
Since $\K$ is infinite, there exists some $z \in \K^\times$ such that  $z^{\la \beta \, |\, \gamma \ra} \neq 1$.
In view of~(\ref{eq:*}), it is now straightforward to check that $z^\gamma$ is  an element of $T$ which
centralizes $U_\alpha$ but not $U_\beta$.
\end{proof}

\subsection{Non-linearity}
\label{ss - NonLinearity}
We can now state the main non-linearity theorem of this section.
Note that it is known that Kac-Moody groups of indefinite type over  infinite fields of arbitrary characteristic do not admit any  \emph{faithful} finite-dimensional linear representation over any field  \cite[Theorem~7.1]{CapracePhD}, but the simplicity for Kac-Moody groups  over infinite fields is still an open question.

\begin{theorem}
\label{thm:NonLinearity} Let $A$ be a generalized Cartan matrix, let $\mathcal{G}_A$ be a Tits  functor of type
$A$ and let $\K$ be a field of characteristic $p>0$. Assume that $A$ is indecomposable, of indefinite type, i.e.
neither  spherical nor affine, that each rank one subgroup of $\mathcal{G}_A(\K)$ is perfect and that
$\mathcal{G}_A(\K)$ is generated by its root subgroups. Then any finite-dimensional linear representation of
$\mathcal{G}_A(\K)$ is trivial.
\end{theorem}

\begin{proof}
We set $G=\mathcal{G}_A(\K)$ and we let $\varphi : G \to \GL_n(\F)$ be  any representation. If $\K$ is finite,
then $\varphi(G)$ is residually finite by Mal'cev's theorem. On the other hand $G$ does not have any finite
quotient by Corollary~\ref{cor:NoFiniteQuotient}. Hence $\varphi(G)$ is trivial in this case.

\vspace{2mm}

We henceforth assume that $\K$ is infinite. Let us denote by $G_p$ the group $\mathcal{G}_A({\bf F}_p)$ (where
${\bf F}_p$ is the prime field of $\K$). We view  $G_p$ as a subgroup of $G$. By Mal'cev's theorem,
$\varphi(G_p)$ is residually finite, so by Corollary~\ref{cor:NoFiniteQuotient} it is finite. Since the root
system $\Phi$ contains nilpotent subsets of arbitrary  large cardinality, the kernel $H$ of $\varphi$ meets
non-trivially the  ${\bf F}_p$-points of $U_\Psi$ for some nilpotent set of roots $\Psi$. In particular, we have
$H \cap U_\Psi \neq \{1\}$, which implies by Proposition~\ref{prop:normal:InfiniteFields} that $H$ contains
$G^\dagger$. Since the field $\K$ is infinite, we have $G^\dagger = [G, G]$ because the rank one subgroups
$X_\alpha = \la U_\alpha \cup U_{-\alpha} \ra$ are perfect. The conclusion follows since $[G, G]=G$ by
hypothesis.
\end{proof}

\section{Homomorphisms to topological groups}
\label{s - homomorphisms}

In this section we study homomorphisms from Kac-Moody groups to locally  compact groups. In the first result, we
collect some basic facts which show that the  only interesting group homomorphisms from finitely generated
Kac-Moody  groups are those with totally disconnected target. However, the main part of this section is devoted
to proving that any nontrivial continuous homomorphism whose domain is the topological completion of a twin
building lattice is a proper map. This is a useful result to be combined with superrigidity.

\subsection{Homomorphisms from simple discrete groups}
\label{ss - from discrete simple} We collect here some basic (and probably well-known) facts about abstract
group homomorphisms from simple discrete to locally compact  groups.

\begin{prop}
\label{prop - homomorphisms}
Let $\Lambda$ be an infinite finitely generated group endowed with the  discrete topology.

\begin{enumerate}
\item[(i)] The group $\Lambda$ is residually finite if, and only if, there exists an injective homomorphism from
$\Lambda$ to a compact group.

\item[(ii)] If $\Lambda$ is simple (resp. simple and Kazhdan), any  group homomorphism from $\Lambda$ to a
compact (resp. amenable) group  is trivial.
\end{enumerate}

We henceforth assume that $\Lambda$ is simple.

\begin{enumerate}
\item[(iii)] There exists no nontrivial group homomorphism from  $\Lambda$ to a Lie group with finitely many
connected components.

\item[(iv)] Let $\varphi : \Lambda\to G$ be a nontrivial group  homomorphism to a locally compact group $G$ and
let $\pi : G \to  G/G^\circ$ be the projection onto the group of connected components. Then $\pi \circ \varphi$
is a continuous, injective, unbounded   homomorphism.

\item[(v)] Let $X$ be a ${\rm CAT}(0)$ or hyperbolic proper metric  space.
Then if $\Lambda$ fixes a point, say $\xi$, in the visual boundary  $\partial_\infty X$, it stabilizes each horosphere centered at $\xi$.
\end{enumerate}
\end{prop}

Point (i) was pointed out to us as a folklore result by N.~Monod.

\begin{proof}
(i). By definition, a residually finite group injects in its profinite  completion, so one direction is clear.
Now let $\Lambda$ admit an injective homomorphism $\varphi:\Lambda\to  K$ into a compact group $K$ and let
$\lambda \in \Lambda-\{1\}$. The regular representation $\rho_K$ of $K$ in $L^2(K)$ is injective; we  will use
its Peter-Weyl decomposition \cite[Theorem~27.40]{HewittRoss}.
The image $\bar\lambda$ of $\lambda$ in some suitable finite-dimensional irreducible submodule, say $V$, is
nontrivial. The projection $\Lambda_V$ of $(\rho_K \circ \varphi)(\Lambda)$ to  ${\rm GL}(V)$ is a finitely
generated linear group containing $\bar\lambda$. By Mal'cev's theorem \cite[Window 7 \S 4 Proposition
8]{LubSeg}, the  group $\Lambda_V$ is residually finite, so it admits a finite quotient  in which $\bar\lambda$
is nontrivial: this is a finite quotient of  $\Lambda$ in which the image of the arbitrary nontrivial element
$\lambda$ is nontrivial.

\vspace{2mm}

(ii). The case when $\Lambda$ is simple follows immediately from (i).
We assume that $\Lambda$ is both simple and Kazhdan.
Let $\varphi:\Lambda\to P$ be a homomorphism to an amenable group $P$.
The closure $\overline{\varphi(\Lambda)}$ is a Kazhdan group because so  is $\Lambda$
\cite[Proposition 7.1.6]{Zimmer}
and it is amenable as a closed subgroup of $P$ \cite[2.3.2]{Greenleaf}.
Therefore it is compact \cite[III.3 p.115]{Margulis} and it remains to  apply the first case of this point.

\vspace{2mm}

(iii). Let $\varphi:\Lambda\to G$ be a homomorphism to a Lie group with  finitely many connected components. By
simplicity, the group $\Lambda$ has no finite index subgroup, so we are  reduced to the case when $G$ is
connected. We compose this map with the adjoint representation of $G$,  whose kernel is the center $Z(G)$
\cite[III.6.4 Corollaire 4]{BbkLie2-3},
in order to obtain a continuous homomorphism ${\rm Ad} \circ \varphi$  to the general linear group of the Lie
algebra of $G$. This map is not injective since the group $\Lambda$ is  simple and finitely generated, hence
non-linear. Therefore $({\rm Ad} \circ \varphi)(\Lambda)$ is trivial.  Finally, again by simplicity, we
successively obtain $\varphi(\Lambda)<Z(G)$ and  $\varphi(\Lambda)=\{1\}$.

\vspace{2mm}

(iv). By simplicity of $\Lambda$, the map $\varphi$ is injective since  it is not trivial. Moreover the kernel
of $\pi \circ \varphi$ is equal to $\{1 \}$ or $\Lambda$. We have to  exclude the case when ${\rm
Ker}(\pi\circ\varphi)=\Lambda$. Let us assume the contrary, i.e.  $\varphi(\Lambda)<G^\circ$, in order to obtain
a contradiction. It follows from \cite[4.6]{MZ} that there exists a  compact normal subgroup $K \triangleleft
G^\circ$ such that $G^\circ/K$ is a connected Lie group. Let us  consider the composed map $\Lambda \,\,
{\buildrel \varphi \over \longrightarrow}\,\, G^\circ  \,\, {\buildrel p \over \longrightarrow} \,\, G^\circ/K$
where $p :  G^\circ \to G^\circ/K$ denotes the canonical projection. By  (iii) we have $(p\circ
\varphi)(\Lambda)=\{1\}$ so $\varphi(\Lambda)<K$. It remains to apply  (ii) to obtain the desired contradiction.
The unboundedness of $\pi \circ \varphi$ follows from (ii) as well.

\vspace{2mm}

(v). For each $y \in X$, we denote by $\beta_{\xi,y}$ the Busemann  function $\beta_{\xi,y} : X \to {\bf R}$
centered at $\xi$ and such  that $\beta_{\xi,y}(y)=0$ \cite[II.8.17 and III.H.3]{BH99}. We pick $x  \in  X$ and
define the function $\varphi_{\xi,x} :  \Lambda \to {\bf R}$ by setting $\varphi_{\xi,x}(g)=\beta_{\xi,x}(g.x)$.
Then for $g$, $h  \in  \Lambda$, we compute $\varphi_{\xi,x}(gh)-\varphi_{\xi,x}(h)$, i.e.
$\beta_{\xi,x}(gh.x)-\beta_{\xi,x}(h.x)$ by definition. By equivariance, this is
$\beta_{\xi,x}(gh.x)-\beta_{g.\xi,g.x}(gh.x)$,  that is also $\beta_{\xi,x}(gh.x)-\beta_{\xi,g.x}(gh.x)$ because
$\xi$  is  fixed under the $\Lambda$-action. But the latter quantity is also $\beta_{\xi,x}(g.x)$, i.e.
$\varphi_{\xi,x}(g)$, by the cocycle property of Busemann functions. In other words, the function
$\varphi_{\xi,x}$ is a group homomorphism  from $\Lambda$ to $({\bf R},+)$. By simplicity of $\Lambda$, it is
trivial, from which we deduce that  $x$ and $g.x$ are on the same horosphere centered at $\xi$ for any $g \in
\Lambda$ and any $x \in X$.
\end{proof}

Note that if in (v) we replace $\Lambda$ by a topologically simple  group acting continuously on $X$, the same
conclusion holds (the  argument is the same: the above map $\varphi_{\xi,x}$ is a continuous  group
homomorphism).

\subsection{Diverging sequences in Coxeter groups}
\label{ss - DV in Coxeter} In the present subsection, we consider a Coxeter system $(W,S)$ and the associated
Davis complex $\mathscr{A}$.

\vspace{2mm}

We say that a sequence $(w_n)_{n \geq 0}$ of elements of $W$ \textit{diverges} if $\displaystyle \lim_{n \to
+\infty} \ell(w_n) = +\infty$.

\begin{lem}\label{lem:DivergingCoxeter}
Let $(w_n)_{n \geq 0}$ be a diverging sequence in $W$. Given any $x \in  \mathscr{A}$, there exists a root
half-space $\alpha  \in  \Phi(\mathscr{A})$ and a subsequence  $(w_{n_k})_{k \geq 0}$ such that $\displaystyle
\lim_{k \to +\infty} d(x, w_{n_k}.\alpha) = +\infty$.
\end{lem}

\begin{proof}
The sequence $(w_n)_{n\geq 0}$ diverges if and only if so does  $(w_n\inv)_{n\geq 0}$. Therefore, it suffices to
find a root $\alpha  \in  \Phi$ and a  subsequence $(w_{n_k})_{k \geq 0}$ such that $\displaystyle \lim_{k \to
+\infty} d(w_{n_k}.x, \alpha) = +\infty$. We set $x_n = w_n.x$. Since $(w_n)_{n\geq 0}$ diverges, we have
$\displaystyle \lim_{n \to  +\infty} d(x, x_n) = +\infty$. Therefore $(x_n)_{n\geq 0}$ possesses a subsequence
$(x_{n_k})_{k\geq  0}$ which converges to a point $\xi$ of the visual boundary  $\partial_\infty \mathscr{A}$.

\vspace{2mm}

Let $\rho : [0, +\infty) \to \mathscr{A}$ be the geodesic ray such that  $\rho(0)=x$ and $\rho(+\infty)=\xi$.
Since $\rho$ is unbounded and since chambers are compact, it follows  that $\rho$ meets infinitely many walls of the Davis complex  $\mathscr{A}$.
On the other hand, the ray $\rho$ is contained in finitely many walls,  otherwise its pointwise stabilizer would be infinite, contradicting the  fact that $W$ acts properly discontinuously on $\mathscr{A}$.
Therefore, there exists a wall $\partial\alpha$ which meets $\rho$ and  such that $\rho$ is not contained in it.
This wall determines two roots, one of which containing no subray of  $\rho$.
We let $\alpha$ be that root: $\alpha \cap \rho$ is a bounded  (nonempty) segment.

\vspace{2mm}

Since $\alpha$ is a closed convex subset of $\mathscr{A}$, the map $d(  \cdot, \alpha) : \mathscr{A} \to \R_+$ is convex \cite[II.2.5]{BH99}  and so is $f:\R_+ \to \R_+ : t \mapsto d(\rho(t), \alpha)$.
Therefore, if $f$ is bounded, it is constant and since $\rho$ meets  $\alpha$ we have $f(t)=0$ for all $t$.
This is excluded because by construction the ray $\rho$ is not  contained in $\alpha$.
Thus $f$ is an unbounded convex function and we deduce that  $\displaystyle \lim_{t \to +\infty} d(\rho(t), \alpha) = +\infty$.
Finally, since $\displaystyle \lim_{t \to \infty} \rho(t) =  \lim_{k  \to \infty} x_{n_k} = \xi$, it follows that
$\{x_{n_k} : k \geq 0\}$ is at finite Hausdorff distance from  $\rho\bigl( [0,+\infty) \bigr)$.
Therefore, we obtain $\displaystyle \lim_{k \to +\infty} d(x_{n_k},  \alpha) = +\infty$ as desired.
\end{proof}

\subsection{Properness of continuous homomorphisms}
\label{ss - properness} We now settle our main properness result for continuous group homomorphisms from
topological completions of twin building lattices.

\begin{theorem}
\label{th:proper} Let $\Lambda$ be a group endowed with a twin root datum $\{U_\alpha\}_{\alpha \in \Phi}$ of
type $(W,S)$ with finite root groups, and such that $\Lambda = \la U_\alpha \; | \; \alpha \in \Phi \ra$. Assume
that $W$ is infinite, that $(W,S)$ is irreducible,  that all root groups are solvable and that property (FPRS)
of Sect.~\ref{sect:FPRS} holds. Let $\overline \Lambda_+$ be the positive topological completion of $\Lambda$
and let $\varphi : \overline\Lambda_+ \to G$ be a continuous nontrivial homomorphism to a locally compact second
countable group $G$. Then $\varphi$ is proper.
\end{theorem}

\begin{remark}
In view of the proof below, we can consider Proposition~\ref{prop:FPRS} and Lemma~\ref{lem:DivergingCoxeter} as
substitutes for results such as the contracting or expanding properties of torus  actions on root groups in the
classical algebraic group case \cite[Lemma 5.3]{BM96}.
\end{remark}

\begin{proof}
Let us assume that $\varphi$ is not proper in order to obtain a  contradiction. There exists a sequence
$(\gamma_j)_{j \geq 0}$ eventually leaving  every compact subset of $\overline\Lambda_+$ and such that
$\displaystyle \lim_{j \to +\infty} \varphi(\gamma_j)$ exists in $G$. Let us recall that we can view the Weyl
group $W$ as the quotient  $\widehat{N_\mathscr{A}}/\Omega_\mathscr{A}$ where $\widehat{N_\mathscr{A}}={\rm
Stab}_{\overline\Lambda_+}(\mathscr{A}_+)$  and $\Omega_\mathscr{A}={\rm
Fix}_{\overline\Lambda_+}(\mathscr{A}_+)$. We also have a Bruhat decomposition:
$\overline\Lambda=\bigsqcup_{w\in W} B_+wB_+$, where $B_+$ is the Iwahori subgroup $B_+={\rm
Fix}_{\overline\Lambda_+}(c_+)$. We use it to write $\gamma_j = k_j.n_j.k'_j$ with $k_j$, $k'_j \in   B_+$ and
$n_j \in  \widehat{N_\mathscr{A}}$. Up to passing to a subsequence, we may assume that $(k_j)_{j \geq 0}$  and
$(k'_j)_{j \geq 0}$ are both converging in the compact open  subgroup $B_+$. We set $w_j =
n_j\Omega_\mathscr{A}$. The hypothesis on $(\gamma_j)_{j \geq 0}$ implies that $(w_j)_{j \geq  0}$ is a
diverging sequence in $W$ and that $\displaystyle \lim_{j \to +\infty} \varphi(n_j)$ exists in $G$. We denote
this limit by $g$. In view of Lemma~\ref{lem:DivergingCoxeter}, up to passing to a subsequence, there exists a
root $\alpha \in \Phi(\mathscr{A}_+)$ such that $\displaystyle \lim_{j \to  +\infty} d(c_+, w_j.\alpha) =
+\infty$. Let $u  \in U_{-\alpha}-\{1\}$. Recall that $n_j. u. n_j\inv \in U_{w_j.(-\alpha)}$ for all $j$.
Therefore by (FPRS) we have: $\displaystyle  \lim_{j \to +\infty} n_j.u .n_j\inv = 1$. Applying $\varphi$ we
obtain: $\displaystyle 1 = \lim_{j \to +\infty}  \varphi(n_j.u .n_j\inv) = g.\varphi(u).g\inv$. Thus we have $u
\in \Ker(\varphi)$. By Proposition~\ref{prop:TopoSimplicity}(i), this implies that $\overline \Lambda_+^\dagger
< \Ker(\varphi)$. Since $\Lambda = \Lambda^\dagger$ by assumption, it follows that $\varphi$ is trivial,
providing the desired  contradiction: $\varphi$ is proper.
\end{proof}

\section{Superrigidity}
\label{s - superrigidity}

In this section, we show that recent superrigidity theorems can be  applied to twin building lattices. They
concern actions on ${\rm CAT}(0)$-spaces. We also derive some consequences: non-linearity of irreducible
cocompact lattices in some Kac-Moody groups, homomorphisms of twin building lattices with Kac-Moody targets,
restrictions for actions on negatively curved complete metric spaces.

\subsection{Actions on ${\rm CAT}(0)$-spaces}
\label{ss - CAT(0) superrigidity} The possibility to apply \cite{Monod06} to irreducible cocompact lattices
enables us to prove the following. Note that the existence of irreducible cocompact lattices in this context is
an open problem.

\begin{prop}
\label{prop - non linearity cocompact} Let $\Lambda$ be a twin building lattice generated by its root groups
$\{U_\alpha\}_{\alpha \in \Phi}$, with associated twinned buildings $\mathscr{B}_\pm$. We assume that the Weyl
group $W$ is infinite, irreducible and non-affine, that root groups are nilpotent and that $Z(\Lambda)=1$. Let
$\Gamma$ be an irreducible cocompact lattice in $\overline\Lambda_+  \times \overline\Lambda_-$. Then any linear
image of $\Gamma$ is finite.
\end{prop}

\begin{remark}
The important assumption here is irreducibility. Indeed, the twin buildings associated to some Kac-Moody groups
are  right-angled Fuchsian. For each of the two buildings $\mathscr{B}_\pm$, the completion
$\overline\Lambda_\pm$ contains a cocompact lattice isomorphic to a  convex cocompact subgroup of the isometry
group of a well-chosen real  hyperbolic space \cite[4.B]{RR06}. Taking the product of two such lattices, we
obtain a cocompact lattice  in $\overline\Lambda_- \times \overline\Lambda_+$ which is linear over  the real
numbers.
\end{remark}

The arguments are classical, so we only sketch the proof.

\begin{proof}
By Proposition~\ref{prop:TopoSimplicity}(i), the groups $\overline \Lambda_\pm$ are topologically simple.
By~\cite[Corollary 1.4]{BSh06}, it follows that  $\Gamma$ is just infinite, i.e. all its proper quotients are
finite. Hence any group homomorphism from $\Gamma$ with infinite image is injective; the same holds for any
finite index subgroup of $\Gamma$. Let ${\bf F}$ be a field with algebraic closure $\overline {\bf F}$ and let
$n \geq 2$ be an integer such that there is an injective group homomorphism $\eta : \Gamma \to {\rm GL}_n({\bf
F})$. We must obtain a contradiction. Let $H$ be the Zariski closure of $\eta(\Gamma)$ in ${\rm GL}_n(\overline
{\bf F})$. We denote by $\Gamma^\circ$ the preimage by $\eta$ of the identity  component $H^\circ$. It is a
finite index normal subgroup of $\Gamma$, so as a lattice in $\overline\Lambda_- \times \overline\Lambda_+$ it
is still irreducible  because $\overline\Lambda_\pm$ is topologically simple. We denote by $R(H^\circ)$ the
radical of $H^\circ$ and by $\pi : H^\circ \to H^\circ/R(H^\circ)$ the natural projection. Then $\pi\circ\eta$
is still injective since otherwise, by the normal  subgroup property for $\Gamma^\circ$, the group $\Gamma$
would be virtually solvable, hence amenable, while $\overline\Lambda_- \times \overline\Lambda_+$ is not. We
thus obtain a semisimple group $G$ over $\overline{\bf F}$ and an  injective group homomorphism $\varphi :
\Gamma^\circ \to G$ with Zariski dense image. We choose an algebraic group embedding in some general linear
group: $G  < {\rm GL}_r$. Being cocompact, the lattice $\Gamma^\circ$ is finitely generated. Taking the matrix
coefficients of the elements of some finite symmetric  generating system implies that $\varphi(\Gamma^\circ)$
lies in ${\rm GL}_r({\bf E})$ for some finitely generated field ${\bf E}$. The group $\Gamma^\circ$ is finitely
generated, linear and non-amenable  so by Tits' alternative it contains a non-abelian free group
\cite{TitsFree}. We can find an element with one eigenvalue, say $\lambda$, of infinite multiplicative order, so
there is a local field $\hat{K}$ with absolute value $\mid\! \cdot \!\mid$ and a field extension $\sigma : {\bf
E} \to  \hat{K}$ such that $\mid\! \sigma(\lambda) \!\mid \, \neq 1$. In particular, the subgroup
$\varphi(\Gamma^\circ)$ is unbounded in $G(\hat{K})$. The map $\varphi : \Gamma^\circ \to G(\hat{K})$ satisfies
the two  conditions required to apply Monod's superrigidity \cite[Corollary 4  and Lemma 59]{Monod06}: the
homomorphism $\varphi$ extends to a continuous homomorphism $\tilde\varphi : \overline\Lambda_- \times
\overline\Lambda_+ \to G(\hat{K})$. By topological simplicity of $\overline\Lambda_\pm$, we obtain an  injective
homomorphism $\Lambda \to G(\hat{K})$, which is impossible since $\Lambda$ is infinite, virtually simple and
finitely generated, hence non-linear.
\end{proof}

\subsection{Uniform $p$--integrability}
\label{ss - uniform integrability} We now check an integrability condition (a substitute for  cocompactness)
required by \cite{GKM06} for a lattice  to enjoy  superrigidity properties. Let us recall the context and the
notation of this reference: $G$ is a  locally compact group, $\Gamma$ is a lattice in $G$. We assume that
$\Gamma$ contains a finite generating subset $\Sigma$  and denote by $\mid\! \cdot \!\mid_\Sigma$ the length
function on $\Gamma$ with respect to it. Following \cite[Sect. 5]{GKM06}, each time we have a right fundamental
domain $\Omega$ for $\Gamma$, we define the function $\chi_{{}_\Omega}  : G \to \Gamma$ by $g  \in
\chi_{{}_\Omega}(g)\Omega$. For each real number $p>1$, we say that $\Gamma$ is {\it  $p$--integrable}~if there
is a right fundamental domain $\Omega$ such  that for any $c  \in  G$, we have: $\displaystyle \int_\Omega
(\mid\! \chi_{{}_\Omega}(gc) \mid_\Sigma)^p  \, {\rm d}g \, < + \infty$. The main result of \cite{Rem05} is that
Kac-Moody lattices are  $p$--integrable for any $p>1$. This amounts to saying that the function $g \mapsto \,\,
\mid\!  \chi_{{}_\Omega}(gc) \!\mid_\Sigma$ belongs to $L^p(\Omega)$ for any $c   \in  G$. We are interested in
a stronger property. We denote by $\parallel\! \cdot \!\parallel_{\Omega,p}$ the $L^p$--norm  of measurable
functions on $\Omega$.

\begin{defi}
\label{def - uniformly integrable} Given $p  \in  [1;\infty)$, the lattice $\Gamma$ in $G$ is called  {\rm
uniformly $p$--integrable}~if there is a right fundamental domain  $\Omega$ as above such that for any compact
subset $C$ in $G$ we have:

\vspace{2mm}

\centerline{
$\displaystyle \sup_{c \in C} \,\, \int_\Omega (\mid\!  \chi_{{}_\Omega}(gc) \mid_\Sigma)^p \, {\rm d}g  \, < +\infty$,}

\vspace{2mm}

i.e. the real valued function $\varphi_{{}_{\Omega,p}} : c  \, \mapsto \, \parallel\! \chi_{{}_\Omega}(\cdot \, c)  \!\parallel_{\Omega,p}$ is bounded on compact subsets of $G$.
\end{defi}

The relation with Y. Shalom's condition \cite{Shalom} is the following. Given a left fundamental domain $D$ for
the inclusion $\Gamma < G$, we  can define a map $\alpha_D : G \times D \to \Gamma$ by setting
$\alpha_D(g,d)=\gamma$  if, and only if, we have $gd\gamma  \in  D$. Up to translating it, we may -- and shall
-- assume that $1_G$ belongs to $D$. The set $\Omega=D^{-1}$ is a right fundamental domain and the following
equivalences hold:

\vspace{2mm}

\centerline{ $\alpha_D(g,1_G) = \gamma \,\,\Longleftrightarrow\,\, g\gamma  \in  D \,\,\Longleftrightarrow\,\,
\gamma^{-1}g^{-1} \in  \Omega \,\,\Longleftrightarrow\,\, g^{-1}  \in  \gamma \Omega \,\,\Longleftrightarrow\,\,
\gamma = \chi_{{}_\Omega}(g^{-1})$. }

\vspace{2mm}

In other words, we have: $\chi_{{}_\Omega}(g)=\alpha_D(g^{-1},1_G)$ and  $2$--integrability amounts to Y.
Shalom's condition (1.5) in \cite[1.II  p.14]{Shalom}.

\vspace{2mm}

Let us now turn to the specific case of twin building lattices. For the rest of the section, we let $\Gamma
=\Lambda$ be a twin building lattice with twin root datum $\{U_\alpha\}_{\alpha \in \Phi}$. We also let
$G=\overline\Lambda_- \times \overline\Lambda_+$. In \cite[1.2]{Rem05} a left fundamental domain $D$ is defined
by means of refined Tits systems arguments. It is a union $D=\bigsqcup_{w \in W} D_w$, indexed by the Weyl group
$W$, of compact open subsets $D_w$ in $G$ and we have $1_G  \in D_{1_W}$. For the rest of the subsection,
$\Sigma$ denotes the finite symmetric generating subset used in \cite[Definition 1]{RR06}, i.e. the union of the
rank one subgroups $X_\alpha \cdot T =\langle U_\alpha, U_{-\alpha} \rangle \cdot T$ indexed by the simple roots
$\alpha \in \Pi$.

\begin{theorem}
\label{th - p-integrable}  There exists a right fundamental domain $\Omega$ such that for any $p   \in
[1;+\infty)$ the function $\varphi_{{}_{\Omega,p}} : c \, \mapsto \, \parallel\!  \chi_{{}_\Omega}(\cdot \, c)
\!\parallel_{\Omega,p}$ is bounded from  above by a function which is constant on each product of double cosets
modulo the standard Iwahori subgroups in $\overline\Lambda_-$ and $\overline\Lambda_+$. In particular, twin
building lattices are uniformly $p$--integrable for any  $p  \in  [1;+\infty)$.
\end{theorem}

\begin{remark}
The second assertion implies  the first one because the standard Iwahori subgroups $B_\pm$ are open and compact,
so products of double cosets modulo Iwahori subgroups $B_+w_+B_+ \times B_-w_-B_-$ are open and compact (and
disjoint when distinct).
\end{remark}

\begin{proof}
We first note that as for the normal subgroup property (Theorem \ref{th - NSP}), though the results in
\cite{Rem05} are stated for Kac-Moody groups, we can use them in the more general context of the above statement
thanks to Proposition~\ref{prop:TopoCompletions}(vi). Moreover since the groups $\overline\Lambda_\pm$ have
$BN$--pairs, they are unimodular \cite[IV.2.7]{BbkLie4-5-6}. Hence so is $G$. We denote by ${\rm d}g$ a Haar
measure on $G$ and compute  $\varphi_{{}_{\Omega,p}}(c)^p$ for $c  \in  G$. It is:

\vspace{2mm}

\centerline{
$\displaystyle
\int_\Omega \bigl( \mid\! \chi_{{}_\Omega}(gc) \!\mid_\Sigma \bigr)^p  \, {\rm d}g
= \int_{D^{-1}} \bigl(\mid\! \alpha_D(c^{-1}g^{-1},1_G) \!\mid_\Sigma \bigr)^p  \, {\rm d}g
= \int_{D}  \bigl(\mid\! \alpha_D(c^{-1}g,1_G) \!\mid_\Sigma \bigr)^p \, {\rm  d}g$.}

\vspace{2mm}

The first equality follows from the remarks before the statement and  the last equality from the unimodularity
of $G$. Therefore it is enough to check that the map $\displaystyle h \mapsto \int_{D} \ell_\Sigma \bigl(
\alpha_D(hg,1_G)  \bigr)^p \, {\rm d}g$ is bounded from above by a function which is  constant on products of
double cosets modulo the standard Iwahori subgroups $B_-$ and $B_+$, in $\overline\Lambda_-$ and
$\overline\Lambda_+$ respectively. We are now back to objects studied in \cite{Rem05}. An element $h  \in  G$ is
a couple $(h_-,h_+)$ with $h_\pm  \in   \overline\Lambda_\pm$ and we denote by $L_\pm(h)$ the combinatorial
distance (i.e. the length of a minimal gallery) in the building $\mathscr{B}_\pm$ from the standard chamber
$c_\pm$ to the chamber $h_\pm^{-1}.c_\pm$. The function $h \mapsto L_\pm(h)$ is constant on each double coset of
$\overline\Lambda_\pm$ modulo the standard Iwahori subgroup of sign  $\pm$ since $L_\pm(h)$ is nothing but the
length in $W$ of the Weyl group element indexing the double class $B_\pm h_\pm B_\pm$. We set
$L(h)=L_-(h)+L_+(h)$ and introduce the polynomial $Q_h$ defined  in \cite[Lemma 17]{Rem05} by
$Q_h(X)=3X^2+(6L(h)+3)X+(3L(h)^2+3L(h)+1)$. Then by the proof of the main theorem of [loc. cit., p.~39] there
exists a  constant $\mid\! T \!\mid$ such that we have: $\displaystyle \varphi_{{}_{\Omega,p}}(h)^p \,\, \leq
\,\, \mid\! T  \!\mid \cdot \sum_{n \in {\bf N}} {Q_h(n)^p \over q^n}$. We are done since $h \mapsto Q_h$ is
constant on the products of double cosets $B_+w_+B_+ \times B_-w_-B_-$.
\end{proof}

The consequences of the above integrability property are the  possibility to apply the various rigidity
statements of \cite{GKM06}. One of them is explicitly stated in the introduction of the present paper.

\begin{cor}
\label{cor - harmonic superrigidity} Let $\Lambda$ be a twin building lattice. Then theorems {\rm 1.1}, {\rm
1.3} and {\rm 1.4} of {\rm \cite[Introduction]{GKM06}} can be applied to $\Lambda$.
\end{cor}

\begin{remark}
In this subsection, no assumption was made on the type of the Coxeter diagram of the Weyl group.
\end{remark}

\subsection{Homomorphisms of twin building lattices with Kac-Moody targets}
\label{ss - KM target} The purpose of this subsection is to present a concrete application of superrigidity of
twin building lattices, and more specifically of Corollary~\ref{cor - harmonic superrigidity}. Recall that the
main application of superrigidity of lattices in Lie groups is arithmeticity, see e.g. \cite{Margulis} and
\cite{Monod06}. In the context of twin building lattices, and in view of the simplicity
theorem~\ref{thm:SimpleDRJ} and its corollaries, it is rather natural to apply superrigidity to homomorphisms
with non-linear targets. The main result of this section is an example of such an application.

\vspace{2mm}

Let $(G, \{U_\alpha\}_{\alpha \in \Phi})$ be a twin root datum of type $(W,S)$, with finite root groups, such
that the centralizer $Z_G(G^\dagger)$ is trivial, and let $\overline G_+$ be its positive topological
completion. Since $Z_G(G^\dagger)$ is the kernel of $\overline G_+$-action on the associated building
$\mathscr{B}_+$, we may -- and shall -- view $\overline G_+$ as a subgroup of $\Aut(\mathscr{B}_+)$ (see
Proposition~\ref{prop:TopoCompletions}(iii)).

\begin{theorem}\label{thm:KM target}
Let $\Lambda$ be a twin building lattice of irreducible type, whose root groups are solvable, and such that
$\Lambda$ is generated by the root groups and that property (FPRS) of Sect.~\ref{sect:FPRS} holds.  Let $\varphi
: \Lambda \to \overline G_+$ be a homomorphism with dense image. Assume that $(W,S)$ is irreducible,
non-spherical and $2$--spherical (i.e. any $2$--subset of $S$ generates a finite group). Assume also that
$q_{\min} = \min\{|U_\alpha| \; : \; \alpha \in \Phi\} > 3$. Then $\varphi(\Lambda)$ is conjugate to $G^\dagger$
in $\Aut(\mathscr{B}_+)$. In particular $G^\dagger$ is isomorphic to $\Lambda/Z(\Lambda)$.
\end{theorem}
\begin{proof}
Details of the arguments involve some rather delicate considerations pertaining to the theory of twin buildings.
Since a detailed proof would therefore require somewhat lengthy preparations which are too far away from the
main topics of this paper, we only give a sketch.

\vspace{2mm}

Since the $\overline G_+$-action on the associated building $\mathscr{B}_+$ is strongly transitive and since
$\overline{\varphi(\Lambda)}$ is dense in $\overline G_+$, it follows that the $\Lambda$-action on
$\mathscr{B}_+$ induced by $\varphi$ is reduced (in the sense of \cite{Monod06}). By Corollary~\ref{cor -
harmonic superrigidity}, we may therefore apply \cite[Theorem~1.1]{GKM06}, which ensures that the homomorphism
$\varphi$ extends uniquely to a continuous homomorphism $\overline \Lambda_+ \times \overline \Lambda_- \to
\overline G_+$, factoring through $\overline \Lambda_+$ or $\overline \Lambda_-$. Up to exchanging $+$ and $-$,
we assume that $\varphi$ extends to a continuous homomorphism $\overline \Lambda_+ \to \overline G_+$, also
denoted $\varphi$. Since $\varphi$ is proper by Theorem~\ref{th:proper}, it follows that  $\varphi$ is
surjective. Furthermore, by Proposition~\ref{prop:TopoSimplicity}(i) the kernel of $\varphi$ is contained in the
discrete center $Z(\overline \Lambda_+) = Z(\Lambda)$.

\vspace{2mm}

For both $\overline \Lambda_+$ and $\overline G_+$, the maximal compact subgroups are precisely the maximal
(spherical) parahoric subgroups. Moreover, since $(W,S)$ is irreducible and $2$--spherical, a maximal parahoric
subgroup of $\overline G_+$ has a Levi decomposition as semi-direct product $L \ltimes \overline U$, where $L$
is a finite group of Lie type and $\overline U$ is a pro--$p$ group for some prime $p$ depending only on $G$
(see \cite[Theorem~1.B.1(ii)]{Rem04} and Proposition~\ref{prop:pronilpotency}(i)). Furthermore, in the
decomposition $L \ltimes \overline U$, the \og congruence subgroup\fg $\overline U$ is characterized as the
maximal normal pro--$p$ subgroup. Using the fact that maximal parahoric subgroups of $\overline \Lambda_+$ also
admit Levi decompositions, it is not difficult to deduce that the Levi factors in $\overline \Lambda_+$ are also
finite groups of Lie type and, then, that $\varphi$ induces isomorphisms between the Levi factors in $\overline
\Lambda_+$ and in $\overline G_+$. In view of the description of the respective buildings of the latter groups
as coset geometries modulo parahoric subgroups, this in turn implies that $\varphi$ induces an isomorphism
between the building $\mathscr{X}_+$ of $\overline \Lambda_+$ and the building $\mathscr{B}_+$ of $\overline
G_+$; moreover this isomorphism is $\varphi$--equivariant. In particular $\mathscr{X}_+$ is of $2$--spherical
irreducible type with infinite Weyl group. By assumption, the building $\mathscr{X}_+$ (resp. $\mathscr{B}_+$)
admits a twin $\mathscr{X}_-$ (resp. $\mathscr{B}_-$) such that the diagonal action of $\Lambda$ (resp. $G$) on
the product $\mathscr{X}_+ \times \mathscr{X}_-$ (resp. $\mathscr{B}_+ \times \mathscr{B}_-$) preserves the
twinning. By the main result of \cite{MR94}, these twinnings must be isomorphic since $\mathscr{X}_+$ and
$\mathscr{B}_+$ are. More precisely, these twinnings are conjugate under some element of $\Aut(\mathscr{B}_+)$.
Up to conjugating $\varphi(\Lambda)$ by this element of $\Aut(\mathscr{B}_+)$, we may -- and shall -- assume
that $\varphi(\Lambda)$ preserves the twinning between $\mathscr{B}_+$ and $\mathscr{B}_-$.

\vspace{2mm}

Since the isomorphism between $\mathscr{X}_+$ and $\mathscr{B}_+$ is $\varphi$--equivariant, it follows from
standard description of root group actions in Moufang twin buildings (see e.g. \cite{TitsLMS}) that $\varphi$
maps each root group of $\Lambda$ to a root group of $G$, and that every root group of $G$ is reached in this
manner. Since $\Lambda$ is generated by its root groups, we deduce that $\varphi(\Lambda) = G^\dagger$. Finally,
since $G^\dagger$ is center-free, so is $\varphi(\Lambda)$, from which it follows that $\Ker(\varphi) \cap
\Lambda = Z(\Lambda)$.
\end{proof}

\subsection{Actions on ${\rm CAT}(-1)$-spaces}
\label{ss - CAT(-1) superrigidity} Another consequence is the existence of strong restrictions on actions  of
\og higher-rank\fg Kac-Moody lattices on ${\rm CAT}(-1)$-spaces. Of course, in this case we must discuss the
notion of rank which is  relevant to the situation; this is done just after the statement. Since we are dealing
with hyperbolic target spaces, it is most convenient  to  use a superrigidity theorem due to N. Monod and Y.
Shalom \cite{MS04}.

\begin{theorem}
\label{th - action on CAT(-1)} Let $\Lambda$ be a split or almost split adjoint Kac-Moody group over $\F_q$
which is a lattice in the product of the associated buildings $\mathscr{B}_\pm$. We assume that the Weyl group
$W$ is infinite, irreducible and non-affine and that $q\geq 4$. Let $Y$ be a proper ${\rm CAT}(-1)$-space with
cocompact isometry group. We assume that $\Lambda$ acts on $Y$ by isometries and we denote by  $\varphi :
\Lambda \to {\rm Isom}(Y)$ the corresponding homomorphism.

\begin{enumerate}
\item[(i)] If the $\Lambda$-action is nontrivial but has a global  fixed point in the compactification
$\overline Y= Y \sqcup \partial_\infty Y$, then the fixed point is  unique and lies in the visual boundary
$\partial_\infty Y$.

\item[(ii)]We assume that the $\Lambda$-action has no global fixed  point at all. Then there exists a nonempty,
closed, convex, $\Lambda$-stable subset  $Z \subseteq Y$ on which it extends to a continuous homomorphism
$\tilde\varphi : \overline\Lambda_- \times \overline\Lambda_+ \to {\rm  Isom}(Z)$ which factors through
$\overline\Lambda_-$ or  $\overline\Lambda_+$.

\item[(iii)]We assume that the buildings $\mathscr{B}_\pm$ contain flat subspaces of dimension $\geq 2$. Then
the $\Lambda$-action, if not trivial, has a unique global fixed  point in $\overline Y$, which lies in the
visual boundary $\partial_\infty Y$.

\item[(iv)] If the buildings $\mathscr{B}_\pm$ contain flat subspaces of  dimension $\geq 2$ and if $\Lambda$ is
Kazhdan, then the $\Lambda$-action is trivial.
\end{enumerate}
\end{theorem}

\begin{remarks}
1. The assumption that $Y$ has cocompact isometry group in the theorem is necessary. Consider indeed the minimal
adjoint Kac-Moody group $\Lambda = \mathcal{G}_A(\F_q)$ over $\F_q$, where $A$ is the generalized Cartan matrix
of size~$4$ defined by $A_{ii} = 2$ for $i = 1, \dots, 4$ and $A_{ij} = -1$ for $1 \leq i \neq j \leq 4$. Thus
the group $\Lambda$ satisfies all hypotheses of the theorem if $q \geq 4$; furthermore $\Lambda$ has
property~(T) provided $q > 1764^4$ \cite{DJ02}. The specificity is here that the Weyl group $W$ is a Coxeter
group which is a non-uniform lattice of ${\rm SO}(3,1)$.
In fact $W$ acts on the real hyperbolic $3$-space
$\mathbb{H}^3$ with a non-compact simplex as a fundamental domain. It turns out that the whole building
$\mathscr{B}_+$ has a geometric realization $X$ in which apartments are isomorphic to the tiling of
$\mathbb{H}^3$ by the above simplex as fundamental tile. This geometric realization is a locally finite
simplicial complex which admits a global $\mathrm{CAT}(-1)$--metric induced by the metric of the apartments
\cite[Proposition~11.31]{AB06}. Hence, although the building $\mathscr{B}_+$ has $2$--dimensional flats in its
usual geometric realization, the $\mathrm{CAT}(-1)$-space $X$ is endowed with a natural $\Lambda$-action which
has no global fixed point in the visual compactification $\overline X$. The point is of course that chambers are
not compact in $X$ (while they are of course compact in $\mathscr{B}_+$), which implies that $\Isom(X)$ is not
cocompact.

Note that in the above example, the rank~$2$ Levi factors of $\Lambda$ are (virtually) nothing but arithmetic
groups $\SL_3(\F_q[t,t\inv])$. The action of these subgroups on $X$ induced by the $\Lambda$-action has no
global fixed points in $X$, which shows in particular that the assumption that $\Isom(Y)$ has finite critical
exponent is necessary in \cite[Corollary 0.5]{BM96} as well.

\vspace{2mm}

2. The prototype of \og higher-rank versus CAT(-1)\fg result we have in  mind is \cite[Corollary 0.5]{BM96}. In
the latter case the target space is also a CAT(-1)-space without any  required connection with Lie groups, but
the irreducible lattice lies  in a product of algebraic groups. Then the fact that each factor in the associated
product of  symmetric spaces and Bruhat-Tits buildings contains higher-dimensional  flats implies  property (T)
for the lattice. In the Kac-Moody case, the existence of higher-dimensional flats no longer implies property
(T); this explains the distinction between (iii) and (iv).

\vspace{2mm}

3. The notion of flat rank considered in \cite{BRW05} for groups of  building automorphisms is relevant here.
According to \cite{CapraceHaglund} and \cite[Theorem A]{BRW05}, knowing whether the buildings $\mathscr{B}_\pm$
contain higher-dimensional flat subspaces is equivalent to the fact that the groups $\overline\Lambda_\pm$ are
of flat rank $\geq 2$ or, still equivalently, that the Weyl group contains a free abelian subgroup of rank~$\geq
2$. In particular, the flat rank can be explicitly computed from the Dynkin diagram of the twin building lattice
group $\Lambda$.

\end{remarks}

\begin{proof}[Proof of Theorem~\ref{th - action on CAT(-1)}]
(i). We first note that the hypotheses (S0)--(S3+) of Theorem~\ref{thm:SimpleDRJ} are fulfilled. Moreover
$\Lambda = \Lambda^\dagger$ and $Z(\Lambda)=1$ since $\Lambda$ is adjoint. Hence the group $\Lambda$ is simple.
In particular, the non-triviality of $\varphi$ implies that it is injective. Moreover for any $y  \in Y$ the
stabilizer ${\rm Stab}_{{\rm Iso}(Y)}(y)$ is a compact group, so by Proposition \ref{prop - homomorphisms} a
nontrivial $\Lambda$-action on $Y$ cannot have any  fixed point in $Y$. Now let $\xi$ and $\eta$ be two distinct
$\Lambda$-fixed points in the visual boundary $\partial_\infty Y$. Then the unique geodesic $(\eta \xi)$ is
stable and by simplicity of $\Lambda$ the restriction of the $\Lambda$-action on $(\eta \xi)$ has  to be
trivial: this implies the existence of a global fixed point in $Y$, which again is excluded when $\varphi$ is
nontrivial.

\vspace{2mm}

(ii). Let us first show that the closure group  $\overline{\varphi(\Lambda)}$ is non-amenable. Assume the
contrary. Then there is a probability measure $\mu$ which is fixed by this group. Since the $\Lambda$-action has
no global fixed point in $Y$, the group  $\overline{\varphi(\Lambda)}$ is not compact so by the ${\rm CAT}(-1)$
Furstenberg's lemma \cite[Lemma 2.3]{BM96} the support  of $\mu$ contains at most two points. This support is
$\Lambda$-stable, so by simplicity it is pointwise  fixed by $\Lambda$; this is excluded because the
$\Lambda$-action has  no global fixed point in $\partial_\infty Y$.

\vspace{2mm}

We henceforth know that $\overline{\varphi(\Lambda)}$ is non-amenable.
We apply \cite[Theorem 1.3]{MS04}: there exists a compact normal  subgroup
$M \triangleleft \overline{\varphi(\Lambda)}$ such that the induced  homomorphism
$\Lambda \to \overline{\varphi(\Lambda)}/M$ extends to a continuous  homomorphism
$\overline\Lambda_- \times \overline\Lambda_+ \to  \overline{\varphi(\Lambda)}/M$
factoring through $\overline\Lambda_-$ or $\overline\Lambda_+$.
Therefore, choosing a suitable sign we obtain an injective continuous group homomorphism

\vspace{2mm}

\centerline{$\bar\varphi : \overline\Lambda_\pm \to  \overline{\varphi(\Lambda)}/M$.}

\vspace{2mm}

Let us denote by $Z$ the fixed-point set of $M$ in $Y$. The subset $Z$ is nonempty because $M$ is compact, it is
closed  and convex by uniqueness  of geodesic segments. Since $M$ is normal in $\overline{\varphi(\Lambda)}$,
the subset $Z$ is  stable under the $\overline{\varphi(\Lambda)}$-action. The restriction map $r_Z :
\overline{\varphi(\Lambda)} \to {\rm  Isom}(Z)$ factors through the canonical projection $\pi_M :
\overline{\varphi(\Lambda)} \to \overline{\varphi(\Lambda)}/M$  and provides a natural homomorphism $\bar r_Z :
\overline{\varphi(\Lambda)}/M \to {\rm Isom}(Z)$ which we  can compose with $\bar\varphi$ to obtain the desired
homomorphism  $\tilde\varphi = \bar r_Z \circ \bar\varphi$.

\vspace{2mm}

(iii). Let us assume that the $\Lambda$-action has no fixed point in  $\overline Y$ in order to obtain a
contradiction. By applying (ii), we have an injective continuous homomorphism  $\tilde\varphi :
\overline\Lambda_\pm \to {\rm Isom}(Z)$ as above. Moreover Proposition~\ref{prop:FPRS} and
Theorem~\ref{th:proper} imply that $\tilde\varphi$ is proper.

\vspace{2mm}

{\it Root system preliminaries}.---~ By the remarks before the proof, the existence of flats of dimension $\geq
2$ implies the existence of an abelian subgroup isomorphic to ${\bf Z} \times {\bf Z}$ in the Weyl group $W$.
Moreover  it follows from \cite[Theorem~6.8.3]{Kra94} that if $W$ has a subgroup isomorphic to $\Z \times \Z$,
then it has a reflection subgroup isomorphic to $D_\infty \times D_\infty$, where $D_\infty$ is the infinite
dihedral group. Let $\alpha, \alpha', \beta, \beta' \in \Phi$ be roots such that $\tau = r_\alpha r_\beta$ and
$\tau' = r_{\alpha'} r_{\beta'}$ are mutually commuting and both of infinite order. Let $V_+$ (resp. $V_-$) the
group generated by the root groups indexed by roots in $\bigcup_{n \in {\bf Z}} \tau^n.\{ \alpha;-\beta \}$
(resp. by  the opposite roots). Note that each of these two groups is normalized by any element lifting $\tau$
in $N$.

\vspace{2mm}

{\it Reduction to hyperbolic isometries}.---~ Recall that the torus $T$ is finite, hence the fixed-point-set
$Y^T$ of $T$ in $Y$ is nonempty. Thus $Y^T$ is a closed convex subset of $Y$ on which the group $N$ acts since
$N$ normalizes $T$. Recall that the quotient $N/T$ is nothing but the Weyl group $W$. Obviously $T$ acts
trivially on $Y^T$ and, hence, the action of $N$ on $Y^T$ factors through $W$. Therefore, we may -- and shall --
consider that $W$ acts on $Y^T$.

\vspace{2mm}

Let us now pick $n_\tau$ such an element, i.e. such that $n_\tau T=\tau$ in $W=N/T$. By properness of
$\tilde\varphi$, the group generated by $\tilde\varphi(n_\tau)$ is unbounded, therefore the isometries
$\tilde\varphi(n_\tau)^{\pm 1}$ are either both hyperbolic or both parabolic because $n_\tau$ and $n_\tau \inv$
are mutually conjugate by $r_\alpha$. We claim that we obtain the desired contradiction if we manage to prove
that these isometries (as well as $n_{\tau'}$) are hyperbolic. Indeed, $\tau$ together with $\tau'$ generate a
free abelian group of rank~$2$ which acts on $Y^T$. Since the group $N$ is a discrete subgroup of $\overline
\Lambda_+$ (because it acts properly discontinuously on $\mathscr{B}_+$) and since $\tilde\varphi$ is proper, it
follows that $ \la \tau, \tau' \ra$ acts freely on $Y^T$. By the flat torus theorem \cite[Corollary
II.7.2]{BH99}, we deduce that $Y^T$ contains a $2$-flat. This is absurd because $Y$ is $\mathrm{CAT}(-1)$.

\vspace{2mm}

{\it Fixed points of \og unipotent\fg subgroups}.---~ On the one hand, we claim that
$\overline{\tilde\varphi(V_\pm)}$ cannot  stabilize any geodesic line in $Z$. Indeed, any element $g  \in V_\pm$
is torsion so it fixes a point in  $Z$. If $L$ were a $\tilde\varphi(V_\pm)$-stable geodesic line then, using
orthogonal projection, $\tilde\varphi(g)$ would fix a point of $L$. This would imply that the subgroup of index
at most 2 in $\tilde\varphi(V_\pm)$ fixing the extremities of $L$ would in fact fix  the whole line $L$: this is
excluded because, by properness of  $\tilde\varphi$, the groups $\overline{\tilde\varphi(V_\pm)}$ are not
compact. On the other hand, the groups $V_\pm$ are metabelian \cite[3.2 Example 2]{RemNewton} so the closures
$\overline{\tilde\varphi(V_\pm)}$ are amenable groups  \cite[4.1.13]{Zimmer}, hence fix a probability measure on
$\partial_\infty Z$. By \cite[Lemma 2.3]{BM96} the support of such a measure contains at most two points. By the
previous point, the support must consist of one single point at  infinity and the same argument shows that this
point is the unique  $\tilde\varphi(V_\pm)$--fixed point in $\partial_\infty Z$: we call it $\eta_\pm$.

\vspace{2mm}

{\it Images of translations are not parabolic}.---~ Since $n_\tau$ normalizes $V_\pm$ and since
$(\partial_\infty  Z)^{\tilde\varphi(V_\pm)}=\{ \eta_\pm \}$, the isometry  $\tilde\varphi(n_\tau)$ fixes
$\eta_-$ and $\eta_+$. In order to see that $\tilde\varphi(n_\tau)$ is hyperbolic, it suffices  to show that
$\eta_- \neq \eta_+$. Let us assume that $\eta_-$ and $\eta_+$ are the same boundary point,  which we call
$\eta$. We need to obtain a very last contradiction. Let us consider the group $H=\overline{\langle
V_-,V_+\rangle}$, by definition topologically generated by $V_-$ and $V_+$. It is non-amenable because, as a
completion of a group with twin root datum of  type $\tilde A_1$, it admits a proper strongly transitive action
on a semi-homogeneous locally finite tree. Theorem~\ref{th:proper} implies that the maps $\tilde\varphi$ and
$\bar\varphi$ are proper. On the one hand, properness of $\bar\varphi$ implies that the group $\bar\varphi(H)$
is non-amenable and neither is $\pi_M{}^{-1}\bigl(\bar\varphi(H)\bigr)$ since $\pi_M$ is proper and  surjective.
On the other hand, properness of $\tilde\varphi$ implies that $\tilde\varphi(H)$ is a closed subgroup of ${\rm
Isom}(Z)$, which  implies that $r_Z^{-1}\bigl(\tilde\varphi(H)\bigr)$ is a closed subgroup of ${\rm
Stab}_{\overline{\varphi(\Lambda)}}(\eta)$. Moreover we have: $\pi_M{}^{-1}\bigl(\bar\varphi(H)\bigr) <
r_Z^{-1}\bigl(\tilde\varphi(H)\bigr)$, so the non-amenable group $\pi_M{}^{-1}\bigl(\bar\varphi(H)\bigr)$ is a
closed subgroup of ${\rm Stab}_{\overline{\varphi(\Lambda)}}(\eta)$. The contradiction comes from the fact that
${\rm Stab}_{{\rm Isom}(Y)}(\eta)$ is amenable, since ${\rm Isom}(Y)$ acts co-compactly on  $Y$
\cite[Propositions 1.6 and 1.7]{BM96}.

\vspace{2mm}

(iv). The conclusion follows from (iii), the fact that ${\rm Stab}_{{\rm Isom}(Y)}(\xi)$ is amenable for each
$\xi  \in  \partial_\infty Y$ \cite[Propositions 1.6 and 1.7]{BM96} and the fact that $\Lambda$,  being
infinite, simple and Kazhdan, admits no nontrivial homomorphism  to an amenable group (Proposition \ref{prop -
homomorphisms}).
\end{proof}

Note that the assumption that $\Lambda$ is Kazhdan might be superfluous in (iv). Indeed, we used it only to
exclude that possibility that $\overline{\varphi(\Lambda)}$ has a fixed point in $\partial_\infty Y$ and, hence,
is amenable. But this might also follow using only the fact that $\Lambda$ is infinite, finitely generated,
simple and non-amenable as a discrete group: the question is to determine whether such a group admits a
nontrivial homomorphism into an amenable group. Recall from Proposition~\ref{prop - homomorphisms}(ii) that
every homomorphism of an infinite, finitely generated, simple and non-compact group to a compact group is
trivial. This seems to suggest that the answer to the above question might be negative, which would imply in
particular that the conclusion of (iv) holds even without assuming that $\Lambda$ is Kazhdan.

\bibliographystyle{amsalpha}
\bibliography{SimplSuperrigid}

\vspace{1cm}

D\'epartement de Math\'ematiques
\hfill
Institut Camille Jordan

Universit\'e Libre de Bruxelles
\hfill
UMR 5208 CNRS - Lyon 1

CP 216
\hfill
Universit\'e de Lyon 1

Boulevard du Triomphe
\hfill
21 avenue Claude Bernard

B-1050 Bruxelles, Belgique
\hfill
69622 Villeurbanne Cedex, France

{\tt pcaprace@ulb.ac.be} \hfill {\tt remy@math.univ-lyon1.fr}

\addtolength{\parindent}{-1.6pt}

\end{document}

%% file: NoFiniteQuotient.pstex_t
\begin{picture}(0,0)%
\includegraphics{NoFiniteQuotient.pstex}%
\end{picture}%
\setlength{\unitlength}{1973sp}%
\begingroup\makeatletter\ifx\SetFigFont\undefined%
\gdef\SetFigFont#1#2#3#4#5{%
  \reset@font\fontsize{#1}{#2pt}%
  \fontfamily{#3}\fontseries{#4}\fontshape{#5}%
  \selectfont}%
\fi\endgroup%
\begin{picture}(11295,8421)(64,-7744)
\put(5626,-2161){\makebox(0,0)[lb]{\smash{\SetFigFont{12}{14.4}{\rmdefault}{\mddefault}{\updefault}{\color[rgb]{0,0,0}$\eta$}%
}}}
\put(6001,-5161){\makebox(0,0)[lb]{\smash{\SetFigFont{12}{14.4}{\rmdefault}{\mddefault}{\updefault}{\color[rgb]{0,0,0}$\xi$}%
}}}
\put(2026,-2011){\makebox(0,0)[lb]{\smash{\SetFigFont{12}{14.4}{\rmdefault}{\mddefault}{\updefault}{\color[rgb]{0,0,0}$\alpha$}%
}}}
\put(5176,389){\makebox(0,0)[lb]{\smash{\SetFigFont{12}{14.4}{\rmdefault}{\mddefault}{\updefault}{\color[rgb]{0,0,0}$\tau=r_\eta \circ r_\alpha$}%
}}}
\put(8701,-1036){\makebox(0,0)[lb]{\smash{\SetFigFont{12}{14.4}{\rmdefault}{\mddefault}{\updefault}{\color[rgb]{0,0,0}$\beta = (\tau)^h.(-\alpha)$}%
}}}
\put(9601,-4411){\makebox(0,0)[lb]{\smash{\SetFigFont{12}{14.4}{\rmdefault}{\mddefault}{\updefault}{\color[rgb]{0,0,0}$\tau=r_\xi \circ r_\beta$}%
}}}
\put(5701,-7636){\makebox(0,0)[lb]{\smash{\SetFigFont{12}{14.4}{\rmdefault}{\mddefault}{\updefault}{\color[rgb]{0,0,0}$\gamma = (\tau')^h.(-\beta)$}%
}}}
\end{picture}